\documentclass[english]{amsart}
\usepackage[T1]{fontenc}
\usepackage[latin9]{inputenc}
\usepackage{geometry}
\usepackage{color}
\usepackage{babel}
\usepackage{amstext}
\usepackage{amsthm}
\usepackage{amssymb}
\usepackage[unicode=true,pdfusetitle,
 bookmarks=true,bookmarksnumbered=false,bookmarksopen=false,
 breaklinks=false,pdfborder={0 0 0},pdfborderstyle={},backref=false,colorlinks=true]
 {hyperref}
\hypersetup{
 linkcolor=brown, citecolor=blue, pdfstartview={FitH}, hyperfootnotes=false, unicode=true}

\makeatletter
\theoremstyle{plain}
\newtheorem{thm}{\protect\theoremname}
\theoremstyle{definition}
\newtheorem{problem}[thm]{\protect\problemname}
\theoremstyle{plain}
\newtheorem{lem}[thm]{\protect\lemmaname}
\theoremstyle{plain}
\newtheorem{cor}[thm]{\protect\corollaryname}
\theoremstyle{definition}
\newtheorem{defn}[thm]{\protect\definitionname}
\theoremstyle{plain}
\newtheorem{prop}[thm]{\protect\propositionname}
\theoremstyle{remark}
\newtheorem{rem}[thm]{\protect\remarkname}
\theoremstyle{plain}
\newtheorem{fact}[thm]{\protect\factname}
\theoremstyle{definition}
\newtheorem{example}[thm]{\protect\examplename}

\usepackage{tikz-cd}
\usepackage{bbm}

\makeatother

\providecommand{\corollaryname}{Corollary}
\providecommand{\definitionname}{Definition}
\providecommand{\examplename}{Example}
\providecommand{\factname}{Fact}
\providecommand{\lemmaname}{Lemma}
\providecommand{\problemname}{Problem}
\providecommand{\propositionname}{Proposition}
\providecommand{\remarkname}{Remark}
\providecommand{\theoremname}{Theorem}

\begin{document}
\title[Remarks on descriptions of groups]{Remarks and Problems about algorithmic descriptions of groups}
\author{Emmanuel Rauzy}
\begin{abstract}
Motivated by a theorem of Groves and Wilton, we propose the study
of the lattice of numberings of isomorphism classes of marked groups
as a rigorous and comprehensive framework to study global decision
problems for finitely generated groups. We establish the Rice and
Rice-Shapiro Theorems for recursive presentations, and establish similar
results for co-recursive presentations. We give an algorithmic characterization
of finitely presentable groups in terms of semi-decidability of two
decision problems: the word problem and the marked quotient problem,
which we introduce. We explain how this result can be used to define
algorithmic generalizations of finite presentations. Finally, we discuss
how the Adian-Rabin Theorem provides incomplete answers in several
respects.
\end{abstract}

\maketitle

\section{\label{sec:Introduction}Introduction}

The study of decision problems for groups began with Max Dehn who,
in 1911 (\cite{Dehn1911}, see \cite{Dehn1987} for a translation),
formulated the three famous problems which are now associated with
his name: the word problem, the conjugacy problem and the isomorphism
problem. His motivation for introducing these came from topology,
in particular from the study of the fundamental group, which had been
introduced not long before by Poincaré. Because of this, he defined
his problems only for finitely presented groups: his article begins
by stating that his goal is to understand the structure of ``the
general discontinuous groups {[}...{]} given by $n$ generators and
$m$ relations between them''. 

Even though Max Dehn defined the word problem and the conjugacy problem
only for finitely presented groups, this restriction is unnecessary,
because these problems are \emph{local decision problems}: they are
set inside a single group, and concern the elements of this group.
It is not necessary to suppose that a finitely generated group $G$
has a finite presentation, nor any kind of finite description, for
the word problem to make sense in it. 

The situation is different for the isomorphism problem: this is a
\emph{global decision problem}, and for it to make sense, one must
be able to provide finite descriptions of possibly infinite groups. 

\medskip

The goal of this article is to discuss several research directions
that concern global decision problems about finitely generated groups
given by descriptions other than finite presentations. Because finite
presentations of groups provide a very good basis for the study of
global decision problems for groups, this may seem like a gratuitous
endeavor. Our main motivation lies in the following result of Daniel
Groves and Henry Wilton:
\begin{thm}
[\cite{Groves2009}] \label{thm:GrovesWilton-1}There exists an
algorithm that, given as input a presentation for a group $G$ and
a solution to the word problem in $G$, determines whether or not
$G$ is free.
\end{thm}

This result is non-trivial and its proof is not elementary, even though
it answers a basic question. It relies on a good understanding of
the universal theory of free groups (Makanin's algorithm to decide
whether a universal sentence is true in all free groups \cite{Makanin1985})
and of the groups that are models of this theory (the \emph{limit
groups} \cite{Sela2001}), including results of \cite{Bumagin2007,Wilton2008,Dahmani2008}. 

One could go as far as stating that this result is proof of the fact
that the theory of global decision problems for groups should not
be built using finite presentations as the basic type of finite description
of a group, but using instead finite presentations supplemented by
word problem algorithms. 

In any case, this result at least shows that there are interesting
problems about global decision problems for groups outside of finite
presentations. 

\subsection{Undecidability results as non-classification results}

A first remark is that many descriptions other than finite presentations
are already commonly used to study computability with groups: matrices
for linear groups \cite{Detinko2019}, automata \cite{Epstein1992},
L-presentations \cite{Bartholdi2003}, different notions of ``computable
presentations'' for countable groups that need not be finitely generated
\cite{MELNIKOV2014} (these are not presentations in the sense of
generators and relations). 

However, our goal here is to discuss the ways in which it is possible
to obtain a theory of decision problems that serves the same purpose
as the one based on finite presentations, in a very precise sense:
we want to be able to consider that the algorithmic tractability of
a class of groups is a (weak) form of explicit classification of these
groups. This is very much in the spirit of Max Dehn, who introduced
his problems as a step towards developing a comprehensive theory of
finitely presented groups. 

For instance, the isomorphism problem for cyclic groups given as computable
groups is unsolvable \cite{Lockhart1981}. One will of course not
consider that this result proves that there can exist no classification
of the set of cyclic groups. On the other hand, unsolvability of the
isomorphism problem for finitely presented groups shows that there
cannot exist a classification of finitely presented groups similar
to that of finite simple groups. As another example, the fact that
the set of finitely presented residually finite groups is not computably
enumerable (Proposition \ref{prop: RF not CE}) can be considered
as proving that there is no good answer to the question ``what are
the finitely presentable residually finite groups''. 

We summarize this as a problem. While its statement is informal, it
will be useful throughout this article, and we use it only in unambiguous
cases, as above with cyclic groups. 
\begin{problem}
\label{prob: Descriptions that are strong enough}Define descriptions
of groups that contain enough information so that undecidability results
can be seen as non-classification results. 
\end{problem}

For such descriptions, we expect many global decision problems to
be undecidable in general, but to become decidable when restricted
to tame classes of groups (finitely generated abelian groups, finite
groups). 

\subsection{General setting: numberings of marked groups }

Our first step is to provide a framework in which decision problems
for different types of descriptions of groups can be studied rigorously. 

In particular, the theorem of Groves and Wilton quoted above, Theorem
\ref{thm:GrovesWilton-1}, was written here exactly as it appeared
in \cite{Groves2009}, but its statement is not completely unambiguous.
In fact, Groves and Wilton themselves felt this had to be adressed,
and together with Manning, they revisited it in \cite{Groves2012},
proposing a precise notion, that of being \emph{computable modulo
the word problem,} in order to formalize it. We discuss this formalism
and its similarity with Banach-Mazur computability in Section \ref{subsec:Example-of-why-it-is-useful}. 

The general setting we present is \emph{the study of the lattice of
equivalence classes of subnumberings of isomorphism classes of marked
groups. }

The notion of numbering was introduced by Malcev in \cite{Malcev1961}
(translated in \cite{Malcev1971}). A \emph{numbering} of a set $X$
is a partial surjection $\nu:\subseteq\mathbb{N}\rightarrow X$. A
\emph{subnumbering} of $X$ is a numbering of a subset of $X$. A
function $f:X\rightarrow Y$ between subnumbered sets $(X,\nu)$ and
$(Y,\mu)$ is called $(\nu,\mu)$\emph{-computable} when there is
a partial computable map $F:\subseteq\mathbb{N}\rightarrow\mathbb{N}$
(so computable in the original sense of Church and Turing) such that
for all $n$ where $\nu$ is defined, $F$ is also defined and $f(\nu(n))=\mu(F(n))$.

Whenever mathematical objects admit some type of finite descriptions,
it is possible to define an associated numbering, the natural number
mapped to a point $x$ encoding a description of $x$. A natural number
that encodes a description of $x$ is called a \emph{name} of $x$.
The definition of a computable function is then understood as follows:
a function is computable when, given the name of a point, it is possible
to compute a name of its image. 

Subnumberings are considered up to equivalence: two subnumberings
$\nu$ and $\mu$ of $X$ are equivalent when the identity of $X$
is both $(\nu,\mu)$- and $(\mu,\nu)$-computable. The set of equivalence
classes of subnumberings of $X$ is a lattice with meet and join operations
corresponding respectively to the conjunction $\wedge$ and the disjunction
$\vee$: the name of a point for $\nu\wedge\mu$ encodes both a name
of this point for $\nu$ and a name of this point for $\mu$, and
the name of a point for $\nu\vee\mu$ is either a name of this point
for $\nu$, or a name of this point for $\mu$.

In Section \ref{subsec:Numberings-associated-to-REC pres etc}, we
define numberings $\nu_{FP}$, $\nu_{RP}$, $\nu_{co-RP}$ and $\nu_{WP}$,
associated respectively to groups given by finite presentations, recursive
presentations, co-recursive presentations and word problem algorithms. 

Finally, with these, we can give a precise statement for Theorem \ref{thm:GrovesWilton-1}:
freeness is a $\nu_{FP}\wedge\nu_{WP}$-decidable property. 

In fact, all four numberings $\nu_{FP}$, $\nu_{RP}$, $\nu_{co-RP}$
and $\nu_{WP}$ are naturally attached to marked groups. A\emph{ $k$-marked
group} is pair $(G,S)$, where $G$ is a group and $S\in G^{k}$ is
a tuple of elements of $G$ that generate it. 

The fact that the use of marked groups is both natural and useful
is rendered apparent by several results of this article. In particular,
the Rice-Shapiro Theorem for recursive presentations (Theorem \ref{thm:Rice-Shapiro-theorem-for-2})
is interpreted as characterizing the topology of $\nu_{RP}$-semi-decidable
sets as the Scott topology on the lattice of $k$-marked groups (see
Section \ref{sec:Lattice-of-marked}), and the algorithmic characterization
of finitely presented groups (Theorem \ref{thm:Algorithmic-characterization-of-FP-gps})
is an effective version of the fact that finitely presented groups
are the compact elements of this lattice (Lemma \ref{lem: FP iff compact }).
In Section \ref{sec:The-Adian-Rabin-theorem}, we look at the Adian-Rabin
Theorem from the point of view of marked groups, and describe several
ways in which it provides incomplete answers. 

\subsection{Decision problems with respect to local descriptions }

We remark that the distinction between local and global problems is
mirrored by a distinction between descriptions that provide local
information and descriptions that provide global information. 

A description provides local information if it allows to compute with
the elements of the given group. For instance, it could be the solution
to the word problem. 

A description provides global information if it relates the given
group to other groups. For instance, finite presentations of groups
provide global information: a group is given by a presentation $\langle S\,\vert\,R\rangle$
when it is the greatest group in which the relations of $R$ hold.
And thus any group which satisfies the relations of $R$ will be a
quotient of $\langle S\,\vert\,R\rangle$. 

Our first results are aimed at showing that decision problems for
groups given by local information are not going to be able to satisfy
the criterion given in Problem \ref{prob: Descriptions that are strong enough}.
We establish the Rice and Rice-Shapiro theorems for groups given by
recursive presentations. 
\begin{lem}
\label{lem:Two groups rec pres-1}Suppose that $(G,S)$ and $(H,S')$
are two recursively presented marked groups, and that $(H,S')$ is
a strict marked quotient of $(G,S)$. Then no algorithm that takes
as input recursive presentations of either $(G,S)$ or $(H,S')$ can
stop exactly on the presentations that define $(G,S)$. 
\end{lem}

\begin{thm}
[Rice-Shapiro theorem for recursive presentations] \label{thm:Rice-Shapiro-theorem-for-2}If
$P$ is property of marked groups that is semi-decidable from recursive
presentations, then there exists a computably enumerable sequence
of finite presentations, such that a group satisfies $P$ if and only
if it is a marked quotient of a group defined by one of these presentations. 
\end{thm}

We explain that this result shows that the topology generated by semi-decidable
properties of recursively presented groups is the Scott topology on
the set of isomorphism classes of marked groups ordered by morphisms. 
\begin{cor}
[Rice theorem for recursive presentations]There is no non-trivial
decidable marked group property for groups given by recursive presentations. 
\end{cor}

Out of these three results, Lemma \ref{lem:Two groups rec pres-1}
is the most problematic one regarding Problem \ref{prob: Descriptions that are strong enough},
because it can be applied even in restricted settings. 

Mann coined in \cite{Mann1982} the term \emph{co-recursive presentation}:
a group $G$ with a generating family $S$ is co-recursively presented
if there exists an algorithm that enumerates all words of $(S\cup S^{-1})^{*}$
that define non-identity elements of $G$. We establish results about
co-recursive presentations that are similar to these quoted above
about recursive presentations, in particular on the impossibility
to distinguish between a group and a strict quotient of it, and a
Rice theorem which states that there is no decidable property of groups
given by co-recursive presentations. These appear in Section \ref{part: Co-rec presentations}.
We note however that there exist dissymmetries between recursive and
co-recursive presentations (the lattice operations of the set of $k$-marked
groups are computable from recursive presentations, but only the join
is computable from co-recursive presentations (Proposition \ref{prop: Compute Meet and Join from presentations}),
and there is no notion of finite co-presentation).

\subsection{Finite presentations and marked quotient algorithms}

Moving on, it seems that purely local information will not permit
to obtain a theory of decision problems for groups that respects the
criterion described in Problem \ref{prob: Descriptions that are strong enough}.
To obtain descriptions of groups that may have a chance to satisfy
this criterion, we start by analyzing the case of finite presentations.
\begin{defn}
[Marked quotient problem] Let $(G,S)$ be a marked group. The \emph{marked
quotient problem} for $(G,S)$ is the problem of deciding if a marked
group $(H,S')$ given by a recursive presentation is a marked quotient
of $(G,S)$. 
\end{defn}

When a residually recursively presented marked group $(G,S)$ has
a semi-decidable marked quotient problem, an algorithm that semi-decides
this problem will constitute a finite description of $(G,S)$. This
allows to define a numbering $\nu_{MQA}$. We prove: 
\begin{thm}
[Algorithmic characterization of finitely presented groups]\label{thm:Algorithmic-characterization-of-FP-gps}A
marked group $(G,S)$ is finitely presented if and only if it has
semi-decidable word problem and semi-decidable marked quotient problem.
And this result is uniform: it states an equivalence of numberings
\[
\nu_{FP}\equiv\nu_{RP}\wedge\nu_{MQA}.
\]
\end{thm}

Finitely presented groups can thus be characterized in terms of solvability
of decision problems. Note that, as finite presentations of groups
play a central role in the theory of decision problems for groups,
it seems only fair that this theory can in turn account for why this
should be the case. 

Theorem \ref{thm:Algorithmic-characterization-of-FP-gps} is not difficult
to prove, it is obtained as a corollary to the Rice-Shapiro Theorem
for recursive presentations. Its main interest lies in the fact that
it renders explicit the local and global information given by finite
presentations: semi-decidability of the word problem is a purely local
form of information, which is embodied in the numbering $\nu_{RP}$
associated to recursive presentations, whereas semi-decidability of
the marked quotient problem is a purely global form of information. 

And what may be the most important remark of this article is that
there is an imbalance between these two types of informations: semi-decidability
of the word problem a very weak form of local information, whereas
semi-decidability of the marked quotient problem is a very strong
form of global information. 

The study of decision problems for groups given by word problem algorithms
together with finite presentations, as begun by Groves and Wilton
in \cite{Groves2009}, can be seen as a rebalancing of the numbering
of finite presentations: 
\[
\nu_{FP}\equiv\nu_{RP}\wedge\nu_{MQA}
\]
 is replaced by 
\[
\nu_{WP}\wedge\nu_{FP}\equiv\nu_{WP}\wedge\nu_{MQA},
\]
and thus the left hand side numbering is strengthened.

But at the same time, it is also possible to weaken the right hand
side numbering that is associated to marked quotient algorithms. Indeed,
global decision problems for groups are always considered in relation
to their possible restrictions to different classes of groups. When
working in a certain set of groups $\mathcal{C}$, it is very natural
to restrict our attention to the marked quotient problem \emph{relative
to $\mathcal{C}$. }That is to say that we consider the following
problem associated to a marked group $(G,S)$: 
\begin{defn}
[Relative marked quotient problem] Let $(G,S)$ be a marked group
and $\mathcal{C}$ a set of marked groups. The \emph{marked quotient
problem} \emph{for} $(G,S)$ \emph{relative to $\mathcal{C}$} is
the problem of deciding if a marked group $(H,S')$ that belongs to
$\mathcal{C}$ given by a recursive presentation is a marked quotient
of $(G,S)$. 
\end{defn}

There are already several examples in the literature of infinitely
presented groups with marked quotient algorithms relative to certain
classes of groups. For instance, in \cite{HARTUNG2011} and \cite{Bartholdi2008}
it was shown that groups that admit L-presentations have marked quotient
algorithms with respect to finite and nilpotent groups. 

In Section \ref{sec:Marked-quotient-algorithms-Relative}, we investigate
the relation between notions of finite presentation modulo a certain
class of groups and existence of marked quotient algorithms. We show
that the existence of a marked quotient algorithm with respect to
a group variety corresponds exactly to the notion of finite presentation
modulo the laws of the variety (Proposition \ref{prop: Quotient algo for group variety }).
We give the example of the lamplighter group as a group that has a
marked quotient algorithm relative to finite groups while not being
finitely presented as a residually finite group. Finally we ask: 
\begin{problem}
Find a class $\mathcal{C}$ of recursively presented groups such that
all groups in $\mathcal{C}$ admit marked quotient algorithms relative
to $\mathcal{C}$ but such that not all groups in $\mathcal{C}$ are
finitely presented as residually $\mathcal{C}$ groups.
\end{problem}

\section{\label{sec:Lattice-of-marked}Lattice of $k$-marked groups and scott
topology}

Here, we introduce the lattice of $k$-marked groups and the Scott
topology on it. 

\subsection{The lattice of $k$-marked groups}

Fix $k\in\mathbb{N}$. A $k$-\emph{marked group} is a finitely generated
group together with a $k$-tuple of elements that generate it. 

A \emph{morphism of $k$-marked groups} from $(G,S)$ to $(H,S')$
is a group morphism between $G$ and $H$ that maps $S$ to $S'$
respecting the order. Such a morphism is an isomorphism if the underlying
group morphism is a group isomorphism, and marked groups are considered
up to isomorphism. 

Note that there is at most one morphism from a $k$-marked group to
another $k$-marked group, and that if there are morphisms $(G,S)\rightarrow(H,S')$
and $(H,S')\rightarrow(G,S)$, then $(G,S)$ and $(H,S')$ are isomorphic
as marked groups: isomorphism classes of marked groups form a partially
ordered set for the quotient relation.

Let $(\mathcal{G}_{k},\rightarrow)$ be the poset of isomorphism classes
of $k$-marked groups. 

Let $(\mathcal{G},\rightarrow)$ be the poset of isomorphism classes
of marked groups, obtained by taking the disjoint union of the sets
$(\mathcal{G}_{k},\rightarrow)$. 

Note that we consider that two groups marked by families of different
cardinalities are incomparable for $\rightarrow$, contrary to what
is customary (see for instance \cite{Champetier2005}): we do not
identify a $k$-marked group $(G,(s_{1},...,s_{k}))$ with the $k+1$-marked
group $(G,(s_{1},...,s_{k},1_{G}))$. The reason for this is that
the different descriptions of marked groups that we consider all provide
explicitly the number $k$ of generators in a marking. The problem
``is a given marked group a $k$-marked group'' is thus always decidable,
and the set of $k$-marked groups should be clopen in the Scott topology
of $(\mathcal{G},\rightarrow)$. See also Section \ref{subsec: marked group recognition adian rabin }
where the fact that we do not identify $(G,(s_{1},...,s_{k}))$ with
$(G,(s_{1},...,s_{k},1_{G}))$ plays a role. 

The poset $(\mathcal{G}_{k},\rightarrow)$ is in fact a lattice: any
two marked groups admit both a $\sup$ and an $\inf$ with respect
to the order. We define its lattice operations. 
\begin{prop}
Any pair of $k$-marked groups has an infimum in $(\mathcal{G}_{k},\rightarrow)$. 
\end{prop}

\begin{proof}
Let $(G,S)$ and $(H,S)$ be two $k$-marked groups, generated by
families which we identify via the canonical bijection. Consider two
presentations $\pi_{1}=\langle S\vert R_{1}\rangle$ and $\pi_{2}=\langle S\vert R_{2}\rangle$
that define respectively $(G,S)$ and $(H,S)$. Define $(G,S)\wedge(H,S)$
to be the group given by the presentation $\langle S\vert R_{1},R_{2}\rangle$.
It is immediate to check that this group constitutes the infimum of
$\{(G,S),(H,S)\}$ for $\rightarrow$. 
\end{proof}
\begin{prop}
Any pair of $k$-marked groups has a supremum in $(\mathcal{G}_{k},\rightarrow)$. 
\end{prop}

\begin{proof}
Let $(G,S)$ and $(H,S')$ be two $k$-marked groups, and denote $S=\{s_{1},...,s_{k}\}$
and $S'=\{s_{1}',...,s_{k}'\}$. Denote by $(G,S)\vee(H,S')$ the
subgroup of the group $G\times H$ generated by the elements $\{(s_{1},s_{1}'),...,(s_{k},s_{k}')\}$.
It is easy to check that this group constitutes the supremum of $\{(G,S),(H,S)\}$
for $\rightarrow$. 
\end{proof}
Another way to define the meet is by expressing $(G,S)$ and $(H,S')$
as quotients of a free group: if $(\mathbb{F}_{S},S)$ maps onto $(G,S)$
with kernel $N_{1}$ and onto $(H,S')$ with kernel $N_{2}$, then
$(G,S)\vee(H,S')$ is given by $(\mathbb{F}_{S}/(N_{1}\cap N_{2}),S)$. 

In fact, it is easy to extend both constructions above to arbitrary
suprema and infima, simply by making the same constructions as above,
but with infinitely many marked groups. A lattice is \emph{bounded}
if it has a least and a greatest element.

The following thus follows.
\begin{thm}
The ordered set $(\mathcal{G}_{k},\rightarrow)$ is a complete bounded
lattice with minimum the trivial group (equipped with its unique $k$-marking)
and with maximum a rank $k$-free group (which also has a unique $k$-marking).
\end{thm}

\subsection{Scott topology of the lattice of $k$-marked groups}

A subset $A$ of a partially ordered set $(L,\le)$ is called \emph{Scott
open }if it is an upper set, i.e. $\forall x\in L,\forall y\in A,\,y\le x\implies x\in A$,
and if it is \emph{inaccessible by directed joins}, i.e. if $D\subseteq L$
is a directed set with a supremum $\sup D$, and if there exists $y\in A$
with $y\le\sup D$, then there exists $x\in D$ with $y\le x$.

Define the \emph{way below relation }$\ll$ on $L$ by $x\ll y$ if
and only if for every directed subset $D$ of $L$ that has a supremum
$\sup D$, if $y\le\sup D$, then there exists $z$ in $D$ with $x\le z$. 

An element $x$ of $L$ is \emph{compact} if $x\ll x$. The following
lemma is immediate.
\begin{lem}
\label{lem: FP iff compact }A marked group is compact in $(\mathcal{G}_{k},\leftarrow)$
if and only if it is finitely presented.
\end{lem}

\begin{proof}
Suppose that $(G,S)$ is finitely presented. Then if $(G,S)\rightarrow\bigwedge_{n\in\mathbb{N}}(H_{n},S_{n})$,
there must $N$ such that the finitely many relations that define
$(G,S)$ already hold in some $\bigwedge_{n\le N}(H_{n},S_{n})$,
so that $(G,S)\rightarrow\bigwedge_{n\le N}(H_{n},S_{n})$ holds. 

Suppose now that $(G,S)$ is not finitely presentable, and let $\langle S\,\vert\,r_{1},r_{2},...\rangle$
be an infinite presentation of it. Then 
\[
(G,S)=\bigwedge_{n\in\mathbb{N}}\langle S\,\vert\,r_{1},r_{2},...,r_{n}\rangle,
\]
and yet $(G,S)\rightarrow\langle S\,\vert\,r_{1},r_{2},...r_{n}\rangle$
never holds.
\end{proof}
A \emph{discriminating family} \cite{Cornulier2007} for a group $G$
is a subset of $G$ that does not contain the identity element and
which intersects every non-trivial normal subgroup of $G$. A group
is \emph{finitely discriminable} if it has a finite discriminating
family.

It is clear that a finitely discriminable is compact in $(\mathcal{G}_{k},\rightarrow)$.
We ask:
\begin{problem}
Are the compact elements of $(\mathcal{G}_{k},\rightarrow)$ exactly
the finitely discriminable groups? 
\end{problem}

\section{\label{part: Preliminaries}Lattice of numbering types}

In this section, we define numberings, the conjunction and disjunction
operations on numberings. While using this formalism may seem needlessly
technical, we give in Section \ref{subsec:Example-of-why-it-is-useful}
two examples that justify the fact that at least having access to
this formalism is beneficial.

\subsection{Definitions of the lattice of equivalence classes of numberings}
\begin{defn}
Let $X$ be a set. A \emph{numbering} of $X$ is a surjection $\nu$
that maps a subset $A$ of $\mathbb{N}$ onto $X$. A \emph{subnumbering}
of $X$ is a numbering of a subset of $X$. We denote this by: $\nu:\subseteq\mathbb{N}\rightarrow X$. 
\end{defn}

The use of subnumberings reflects the fact that, while two types of
finite description of marked groups may be applicable to different
sets of marked groups, we can still consider that they live in a common
set, the set of subnumberings of the set of marked groups. 

The pair $(X,\nu)$ is a \emph{subnumbered set}. The domain of $\nu$
is a subset of $\mathbb{N}$ denoted by $\text{dom}(\nu)$. 

The image $\nu(\text{dom}(\nu))$ of $\nu$ is called the set of $\nu$\emph{-computable
points} of $X$, and denoted $X_{\nu}$. Given a point $x$ in $X$,
an integer $n$ such that $\nu(n)=x$ is called a \emph{$\nu$-name}
of $x$. 

Let $(X,\nu)$ and $(Y,\mu)$ be subnumbered sets. A function $f:X\rightarrow Y$
is called $(\nu,\mu)$\emph{-computable} if there exists a partial
computable function $F:\subseteq\mathbb{N}\rightarrow\mathbb{N}$
such that for all $n$ in the domain of $\nu$, $f\circ\nu(n)=\mu\circ F(n)$.
That is to say, there exists $F$ computable which renders the following
diagram commutative:

\begin{center}
\begin{tikzcd}
X \arrow[r, "f"]& Y \\ \mathbb{N}\arrow[u, "\nu"]\arrow[r, "F"]& \mathbb{N}\arrow[u, "\mu"]\end{tikzcd}
\end{center}
\begin{defn}
If $\nu$ and $\mu$ are subnumberings of a set $X$, denote by $\nu\le\mu$
the fact that the identity of $X$ is $(\nu,\mu)-$computable. 
\end{defn}

In other words, given the $\nu-$name of a point of $X$, it is possible
to compute a $\mu$-name of this point. This formalizes the fact that
$\nu$-names contain more information on points than $\mu$-names. 

It is easy to see that $\le$ defines a pre-order. We call \emph{subnumbering
types }the equivalence classes induced by this pre-order. 
\begin{prop}
The set of subnumbering types on $X$ equipped with the order $\le$
is a lattice. 
\end{prop}

Let $(n,m)\mapsto\langle n,m\rangle$ designate Cantor's pairing function. 
\begin{proof}
We define the meet and join operations, which correspond to the disjunction
$\vee$ and to the conjunction $\wedge$. Let $\nu$ and $\mu$ be
subnumberings of $X$. 

Define $\nu\vee\mu$ by setting, for any natural number $k$, $\nu\vee\mu(2k)=\nu(k)$
and $\nu\vee\mu(2k+1)=\mu(k)$. The domain of $\nu\vee\mu$ is the
set $\{2k,k\in\text{dom}(\nu)\}\cup\{2k+1,k\in\text{dom}(\mu)\}$.

Define a subnumbering $\nu\wedge\mu$ by the following:
\[
\text{dom}(\nu\wedge\mu)=\{\langle n,m\rangle\in\mathbb{N},n\in\text{dom}(\nu),\,m\in\text{dom}(\mu),\nu(n)=\mu(m)\},
\]
\[
\forall\langle n,m\rangle\in\text{dom}(\nu\wedge\mu),\,\nu\wedge\mu(\langle n,m\rangle)=\nu(n).
\]

In each case, it is straightforward to check that the defined operations
do constitute meet and join operations.
\end{proof}
\begin{rem}
The conjunction operation is the more important one of the two operations.
Indeed, a function is computable with respect to a disjunction $\nu\vee\mu$
if and only if it is computable for both $\nu$ and $\mu$. Thus the
study of $\nu\vee\mu$ amounts to that of both $\nu$ and $\mu$.
On the contrary, while it is true that if a function is either computable
for $\nu$, or computable for $\mu$, it will be $\nu\wedge\mu$ computable,
this is not an equivalence. A striking example of this fact is given
by the Groves-Wilton theorem quoted in the introduction of this article,
Theorem \ref{thm:GrovesWilton-1}. 
\end{rem}

\subsection{Semi-decidable sets, Ershov topology}

Let $(X,\nu)$ be a numbered set. 

A subset $A\subseteq X$ is \emph{$\nu$-semi-decidable }if there
is an algorithm that, given the $\nu$-name of a point $x$ in $X$,
stops if and only if $x$ belongs to $A$. 

A set is \emph{$\nu$-co-semi-decidable }if its complement is $\nu$-semi-decidable.

A set is $\nu$\emph{-decidable} if it is both semi-decidable and
co-semi-decidable. 

The intersection of finitely many $\nu$-semi-decidable sets is again
$\nu$-semi-decidable, and the union of a sequence of uniformly $\nu$-semi-decidable
sets is again $\nu$-semi-decidable. Thus the set of $\nu$-semi-decidable
sets resembles a topology. The actual topology generated by the semi-decidable
sets on a numbered set $(X,\nu)$ is called the \emph{Ershov topology}
\cite{Spreen1998}. We will characterize the Ershov topology of the
numbering of recursive presentations as the Scott topology on the
lattice of marked groups. The Ershov topology associated to finite
presentations is the discrete topology. 

The set $X$ is \emph{computably enumerable for $\nu$} if it is empty
or if there is a $(\text{id}_{\mathbb{N}},\nu)$-computable surjection
$f:\mathbb{N}\rightarrow X$. This is one of several possible definitions
for ``being computably countable''. This means that it is possible
to algorithmically produce a list that contains at least one $\nu$-name
for each element of $X$.

A $\nu$-\emph{computable sequence }in $X$ is a $(\text{id}_{\mathbb{N}},\nu)$-computable
function $f:\mathbb{N}\rightarrow X$. 

\subsection{\label{subsec:Numberings-associated-to-REC pres etc}Numberings associated
to finite presentations, recursive presentations, etc }

Whenever we consider objects given by finite data, it is possible
to encode this data with natural numbers, and to define the associated
numbering. In most cases, there is little to no benefit in defining
precisely these numberings. Because we are considering numberings
up to equivalence, the choice of a way of encoding the considered
objects is inconsequential. 

In practical implementation of algorithms, the choice of an encoding
for the considered objects becomes important with respect to time
complexity issues. While this is an interesting problem, this is not
a problem that will be solved by precisely defining how to encode
the considered objects with natural numbers from an abstract standpoint,
without using an explicit programming language. This is another reason
to claim that defining precisely numberings is not useful. 

In any case, we still define here precisely the numberings that are
used in this paper, simply to show that it is possible. 

We use a pairing function $(n,m)\mapsto\langle n,m\rangle$, which
is a computable bijection $\mathbb{N}^{2}\rightarrow\mathbb{N}$ with
a computable inverse. Define a numbering $\Delta$ of finite subsets
of $\mathbb{N}$ by the following: 
\[
\Delta(\langle n,m\rangle)=\{u_{1},...,u_{n}\},
\]
where $m=\langle u_{1},\langle u_{2},...,\langle u_{n-1},u_{n}\rangle...\rangle$. 

Fix also a countably infinite alphabet $\mathcal{A}=\{x_{1},x_{2},...\}$.
For each $n$, define $\mathcal{A}_{n}=\{x_{1},...,x_{n}\}$, let
$\mathcal{A}_{n}^{-1}=\{x_{1}^{-1},...,x_{n}^{-1}\}$ be formal inverses
to the elements of $\mathcal{A}_{n}$. Let $(\mathcal{A}_{n}\cup\mathcal{A}_{n}^{-1})^{*}$
be the set of words on the alphabet $\mathcal{A}_{n}\cup\mathcal{A}_{n}^{-1}$.
Finally, let $\theta_{n}$ be a natural bijection $\theta_{n}:\mathbb{N}\rightarrow(\mathcal{A}_{n}\cup\mathcal{A}_{n}^{-1})^{*}$
(for instance ordering $(\mathcal{A}_{n}\cup\mathcal{A}_{n}^{-1})^{*}$
by the shortlex order). 

Composing $\theta_{n}$ and $\Delta$ yields a numbering of finite
subsets of $(\mathcal{A}_{n}\cup\mathcal{A}_{n}^{-1})^{*}$.
\begin{defn}
The numbering $\nu_{FP}$ of finitely presented marked groups is defined
by: 
\[
\nu_{FP}(\langle n,m\rangle)=\langle\mathcal{A}_{n}\,\vert\,\theta_{n}\circ\Delta(m)\rangle
\]
\end{defn}

To define the numbering associated to recursive presentations, we
use the canonical numbering of all partial computable functions, as
first defined by Church Turing and Kleene: $\varphi_{0}$, $\varphi_{1}$,
$\varphi_{2}$... is an enumeration of the set $\mathcal{PR}$ of
partial computable functions. The standard numbering of computably
enumerable subset of $\mathbb{N}$ is given by $W_{i}=\text{dom}(\varphi_{i})$. 
\begin{defn}
We define the numbering $\nu_{RP}$ associated to recursive presentations
by $\nu_{RP}(\langle n,i\rangle)=\langle\mathcal{A}_{n}\,\vert\,\theta_{n}(W_{i})\rangle$. 
\end{defn}

We will similarly define the numbering associated to co-recursive
presentations. Note that, contrary to what the name indicates, a co-recursive
presentation is not a group presentation, in the sense of a set of
generators and of relators. A \emph{co-recursive presentation} \cite{Mann1982}
is an algorithmic enumeration of \emph{all} the words that define
non-identity elements in a given marked group, i.e. an enumeration
of the co-word problem in this marked group. 

A set of words $R\subseteq(\mathcal{A}_{n}\cup\mathcal{A}_{n}^{-1})^{*}$
is the co-word problem of a group if and only if the complement $R^{c}$
of $R$ in $(\mathcal{A}_{n}\cup\mathcal{A}_{n}^{-1})^{*}$ is a normal
subgroup in the free group over $\mathcal{A}_{n}$. 
\begin{defn}
We define the numbering $\nu_{co-RP}$ associated to co-recursive
presentations by:
\[
\text{dom}(\nu_{co-RP})=\{\langle n,i\rangle\in\mathbb{N},\theta_{n}(W_{i})\text{ is a co-word-problem}\};
\]
\[
\nu_{co-RP}(\langle n,i\rangle)=\langle\mathcal{A}_{n}\,\vert\,(\theta_{n}(W_{i}))^{c}\rangle.
\]
\end{defn}

Note that the domain of this numbering is not a computably enumerable
subset of $\mathbb{N}$, this cannot be avoided. 

Recall from the introduction that a marked quotient algorithm for
a marked group $(G,S)$ is an algorithm that stops exactly on the
codes for recursive presentations of marked quotients of $(G,S)$.
We define the associated numbering by relying on our previous definition
of the numbering associated to recursive presentations.
\begin{defn}
\label{def:Marqued quotient algorithm}We define a numbering $\nu_{MQA}$
of recursively presented groups associated to marked quotient algorithms
by: 
\[
\nu_{MQA}(\langle n,i\rangle)=(G,S)\iff(\nu_{RP}(\langle n,j\rangle)\text{ is a marked quotient of }(G,S)\iff j\in W_{i}).
\]
The domain of $\nu_{MQA}$ is defined to be exactly those $\langle n,i\rangle$
that do define a marked quotient algorithm. 
\end{defn}

We prove in Section \ref{sec:Marked-quotient-algorithms-FP} that
the above definition is correct, that is to say that the marked group
$(G,S)$ is indeed uniquely determined by the above formula. 

Finally, a word problem algorithm is simply a recursive presentation
together with a co-recursive presentation. 
\begin{defn}
The numbering $\nu_{WP}$ associated to word problem algorithms is
given by $\nu_{WP}\equiv\nu_{RP}\wedge\nu_{co-RP}$. 
\end{defn}

\subsection{\label{subsec:Example-of-why-it-is-useful}Examples where the explicit
use of numberings is beneficial}

\subsubsection{First example: recursive presentations }

Our first example comes from a result of Lockhart from \cite{Lockhart1981}.
Here is stated the following: ``There is a r.e. class of recursive
presentations with uniformly solvable word problem for which the properties
of freeness and finiteness are unrecognizable''. The groups in questions
are finite and infinite cyclic groups. 

But note that by Lemma \ref{lem:Two groups rec pres}, the word problem
can never be uniformly solvable from recursive presentations in a
class of marked groups that contain a marked group and a strict quotient
of it. This seems to be contradictory with what is written above.

The reason why there is no contradiction is that in the result of
Lockhart, the ``r.e. class of recursive presentations'' is not a
set of \emph{presentations}, but a certain set of \emph{algorithms}
that produce these presentations. The word problem would not be uniformly
solvable on these groups if we allowed arbitrary algorithms that still
produced the correct presentations.

And in fact, the result of Lockhart concerns groups given by word
problem algorithms, and not groups given by recursive presentations.
Formulating this result in terms of numberings makes its ambiguities
disappear.

\subsubsection{Second example: decidability modulo the word problem and Banch-Mazur
computability }

In \cite{Groves2012}, the authors define a precise notion of \emph{group
property recursive modulo the word problem}, in order have a rigorous
concept that formalizes the idea that a marked group is given by a
finite presentations and a solution to its word problem.

Note that using the formalism of numberings, this amounts simply to
defining a numbering via the conjunction operation: if $\nu_{FP}$
is the numbering associated to finite presentations, and $\nu_{WP}$
the numbering associated to word problem algorithms, this is the study
of the decidable properties of $\nu_{WP}\wedge\nu_{FP}$. 

We want to note here that the formalism introduced in \cite{Groves2012}
is very similar to that of \emph{Banach-Mazur computability}.

A real $x$ is \emph{computable} if there is a computable map $f:\mathbb{N}\rightarrow\mathbb{Q}$
which on input $k$ produces a rational $q$ with $\vert x-q\vert<2^{-k}$.
A \emph{computable sequence of computable reals} is a sequence $(u_{n})\in\mathbb{R}^{\mathbb{N}}$
such that there is a computable map $f$ of two variables which, given
$n$ and $k$, produces a rational approximation of the real $u_{n}$
precise within $2^{-k}$. Let $\mathbb{R}_{c}$ be the set of computable
reals. 
\begin{defn}
A function $f:\mathbb{R}_{c}\rightarrow\mathbb{R}_{c}$ is \emph{Banach-Mazur
computable} if it maps computable sequences of computable reals to
computable sequences of computable reals.
\end{defn}

This notion of computability was introduced by Banach and Mazur right
after Church and Turing formalized computability for functions defined
on the natural numbers \cite{Banach1937}. In modern computable analysis,
Banach-Mazur computability is not considered to be ``the good''
notion of computability on the real numbers, although it still constitutes
an important tool to study other notions of computability. For instance,
Hertling showed in \cite{Hertling2002} that the classical effective
continuity theorem of Kreisel-Lacombe-Schoenfield-Ceitin \cite{Ceitin1967}
does not hold for Banach-Mazur computable functions on the real numbers:
this increases our understanding of the original theorem. 

There are two main problems related to the notion of Banach-Mazur
computability. The first is that it does not permit to define a Cartesian
closed category: a Banach-Mazur computable function is not associated
to a finite description, and so we cannot ask what is ``a Banach-Mazur
computable function defined on the set of Banach-Mazur computable
functions''. The second one is that in practice, whenever proving
that some function $f$ is computable, we always prove more than Banach-Mazur
computability: we provide a \emph{single} method to compute the function
$f$, which is uniform on all computable sequences. And thus working
with Banach-Mazur computability, one ends up systematically proving
strictly more than what one states.

The definition of \emph{recursive modulo the word problem} considered
in \cite{Groves2012} is as follows.
\begin{defn}
[\cite{Groves2012}, Definition 2.4]\label{def: GMW def}A class
of finitely presented groups $\mathcal{C}$ is said to be \emph{recursive
modulo the word problem} if, whenever $\mathcal{D}$ is a set of finite
presentations on which the word problem is uniformly solvable, there
exist two computably enumerable sets of finitely presented groups
$\mathcal{X}$ and $\mathcal{Y}$ such that $\mathcal{C}\cap\mathcal{D}\subseteq\mathcal{X}$,
$\mathcal{D}\setminus\mathcal{C}\subseteq\mathcal{Y}$, $\mathcal{C}\cap\mathcal{D}\cap\mathcal{Y}=\emptyset$
and $(\mathcal{C}\setminus\mathcal{D})\cap\mathcal{X}=\emptyset$.
\end{defn}

The common feature of this and of Banach-Mazur computability is the
common non-effective quantification: a problem is considered solvable
if it is solvable on all sets with a certain property (in the first
case, on computable sequences, here on sets where the word problem
is uniformly solvable from finite presentations), but the solutions
on these different sets may be different, and do not have to depend
uniformly on the considered set. 

The difference between the above and Banach-Mazur computability lies
in the fact that the set $\mathcal{D}$ is not supposed to be computably
enumerable. But in fact, we want to argue that the above definition
would be improved by restricting the quantification on $\mathcal{D}$
to c.e. sets of finitely presented groups. Our argument lies in the
following proposition, which follows from the work of Kharlampovich,
Miasnikov and Sapir \cite{Kharlampovich2017}:
\begin{prop}
\label{prop: RF not CE}The set of finitely presented residually finite
groups is not computably enumerable (for $\nu_{FP}$), and furthermore
it cannot be contained in a $\nu_{FP}$-c.e. set of finitely presented
groups with uniformly solvable word problem.
\end{prop}

We will use the following lemma. Recall that time complexity of the
solution to the word problem inside a group with solvable word problem
is a group invariant up to standard asymptotic equivalence of functions
\cite{Sapir2011}.
\begin{lem}
\label{lem: Diagonal complexity}Let $\mathcal{C}$ be a set of finitely
presented groups in which the word problem is uniformly solvable.
If $\mathcal{C}$ is $\nu_{FP}$-c.e., then there exists a computable
function which is a common asymptotic upper bound to the time complexity
for the word problem for groups in $\mathcal{C}$. 
\begin{proof}
By standard embedding methods (see \cite{Darbinyan2015}), it is possible
to embed the restricted direct product of all the groups in $\mathcal{C}$
in a finitely generated group with solvable word problem. The time
complexity of the word problem in this group constitutes the desired
asymptotic upper bound. 
\end{proof}
\end{lem}

\begin{proof}
[Proof of Proposition \ref{prop: RF not CE}] By \cite{McKinsey1943},
the word problem is uniformly solvable in finitely presented residually
finite groups. By \cite{Kharlampovich2017}, there exist finitely
presented residually finite groups of arbitrarily high time complexity
(i.e. above any given computable function). We can conclude by Lemma
\ref{lem: Diagonal complexity}. 
\end{proof}
Remark that the above lemma would also apply to the set of finitely
presented simple groups if the famous Boone-Higman \cite{Belk2023}
conjecture were true:
\begin{prop}
If the Boone-Higman conjecture holds, then the set of finitely presented
simple groups cannot be contained in a c.e. set of finitely presented
groups with uniformly solvable word problem. 
\end{prop}

Now consider the consequences of Proposition \ref{prop: RF not CE}
on Definition \ref{def: GMW def}: to prove that a set $\mathcal{C}$
is recursive modulo the word problem, applying the definition with
$\mathcal{D}$ being the set of residually finite groups, one needs
to find c.e. sets $\mathcal{X}$ and $\mathcal{Y}$ that, in particular,
must cover $\mathcal{D}$. But then the word problem cannot be uniformly
solvable in $\mathcal{X}\cup\mathcal{Y}$, and thus one ends up working
in sets where the word problem algorithm is not uniformly solvable
to show that something is recursive modulo the word problem. 

If on the contrary we modify Definition \ref{def: GMW def}, and only
allow $\mathcal{D}$ to be a c.e. set, then the resulting new definition
of computable modulo the word problem is exactly Banach-Mazur computability.
And, as we explained above, while this notion is useful for theoretical
purposes, it suffers from the fact that in practice, when proving
something Banach-Mazur computable, one always proves that the considered
function is computable for a more restrictive notion of computability.
For instance, one could rewrite the results of \cite{Groves2012},
replacing \emph{computable modulo the word problem} by \emph{computable
for the numbering $\nu_{WP}\wedge\nu_{FP}$}, and leaving all proof
ideas unchanged. The results thus stated would be strictly stronger. 

\section{\label{part:Rice-Shapiro}Rice-Shapiro for recursive presentations}

Here we establish: 
\begin{thm}
[Rice-Shapiro theorem for recursive presentations] \label{thm:Rice-Shapiro-theorem-for}If
$P$ is property of marked groups that is semi-decidable from recursive
presentations, then there exists a computably enumerable sequence
of finite presentations, such that a group satisfies $P$ if and only
if it is a marked quotient of a group defined by one of these presentations. 
\end{thm}

\begin{cor}
The topology generated by sets that are semi-decidable from recursive
presentations is the Scott topology on the lattice of marked groups. 
\end{cor}

And thus:
\begin{cor}
\label{cor:SD rec pres =00003D> quotient stable}Any property of marked
groups that is semi-decidable from recursive presentations is quotient-stable. 
\end{cor}

Finally, the Rice theorem follows immediately from the Rice-Shapiro
theorem.
\begin{cor}
[Rice theorem for recursive presentations]There is no non-trivial
decidable group property for groups given by recursive presentations.
\end{cor}

\begin{proof}
Let $P$ be a decidable property. Either $P$ or its complement contains
the free group. Then, by Corollary \ref{cor:SD rec pres =00003D> quotient stable},
the one that contains it is in fact the set of all marked groups. 
\end{proof}
The proof of Theorem \ref{thm:Rice-Shapiro-theorem-for} follows closely
the usual one \cite{Rogers1987}, we simply exchange the lattice of
subsets of $\mathbb{N}$ by the lattice of marked groups. It is based
on two intermediate lemmas. 
\begin{lem}
\label{lem:Two groups rec pres}Suppose that $(G,S)$ and $(H,S')$
are two recursively presented marked groups, and that $(H,S')$ is
a strict marked quotient of $(G,S)$. Then no algorithm that takes
as input recursive presentations of either $(G,S)$ or $(H,S')$ can
stop exactly on the presentations that define $(G,S)$. 
\end{lem}

In other words, $\{(G,S)\}$ is not a $\nu_{RP}$-semi-decidable subset
of $\{(G,S),(H,S')\}$. 
\begin{proof}
Consider a presentation $\langle S\,\vert\,r_{1},r_{2},...\rangle$
of $(G,S)$ such that $r_{1},r_{2},...$ is a computable sequence,
and a presentation $\langle S'\,\vert\,q_{1},q_{2},...\rangle$ of
$(H,S')$, $q_{1}$, $q_{2}$, ... also being a computable sequence
of relations. Note that in these presentations we can in fact identify
$S$ and $S'$ without changing the considered marked groups. 

Consider an effective enumeration of all partial computable functions
$\varphi_{0}$, $\varphi_{1}$, $\varphi_{2}$, $\varphi_{3}$,...
We write $\varphi_{n}(k)\uparrow^{t}$ if a run of $\varphi_{n}$
on input $k$ does not end in less than $t$ steps, and $\varphi_{n}(k)\downarrow^{t}$
if it halts in exactly $t$ steps. 

For each $n$, we define a computable enumeration $E_{n}=(s_{1},s_{2},...)$
of relations by the following: $s_{t}=r_{t}$ if $\varphi_{n}(n)\uparrow^{t}$,
and $s_{t}=q_{t-t_{0}}$ if $\varphi_{n}(n)\downarrow^{t_{0}}$. 

By the smn-theorem, the sequence $(E_{n})_{n\in\mathbb{N}}$ defines
a computable sequence of recursive presentations. And $E_{n}$ defines
a presentation of $(H,S')$ if and only if $\varphi_{n}(n)\downarrow$,
and a presentation of $(G,S)$ if and only if $\varphi_{n}(n)\uparrow$.
We conclude by unsolvability of the halting problem. 
\end{proof}
This result implies in particular that if the word problem is uniformly
solvable on a set of recursively presented groups, this set does not
contain a group and a strict quotient of it. This condition is not
superfluous, since recursively presented simple groups have uniformly
solvable word problem (and the set of simple group does not contain
a group and a strict quotient of it). 
\begin{lem}
\label{lem:Second-lemma:-sequence Rice Shapiro}Let $(G,S)$ be a
recursively presented marked group that is not finitely presentable,
and let $(r_{0},r_{1}...)$, $r_{i}\in(S\cup S^{-1})^{*}$, be a computable
sequence of relators for it. Let $(G_{n},S)$ be the marked group
given by the truncated presentation $\langle S\vert r_{0},r_{1},r_{2},...,r_{n}\rangle$. 

Then no algorithm that takes as input recursive presentations of marked
groups in $\{(G,S)\}\cup\{(G_{n},S),n\in\mathbb{N}\}$ can stop exactly
on the presentations that define $(G,S)$. 
\end{lem}

In other words, $\{(G,S)\}$ is not a $\nu_{RP}$-semi-decidable subset
of $\{(G,S)\}\cup\{(G_{n},S),n\in\mathbb{N}\}$. 
\begin{proof}
We use as in the proof of Lemma \ref{lem:Two groups rec pres} an
effective enumeration of all partial computable functions $\varphi_{0}$,
$\varphi_{1}$, $\varphi_{2}$, $\varphi_{3}$,...

For each $n$, we define a computable enumeration $E_{n}=(s_{1},s_{2},...)$
of relations by the following: $s_{t}=r_{t}$ if $\varphi_{n}(n)\uparrow^{t}$,
and $s_{t}=r_{t_{0}}$ if $\varphi_{n}(n)\downarrow^{t_{0}}$. 

Thus $E_{n}$ gives a presentation of $(G,S)$ if and only if $\varphi_{n}(n)\uparrow$,
and a presentation of one of the marked groups $(G_{n},S)$ if and
only if $\varphi_{n}(n)\downarrow$. We conclude by unsolvability
of the halting problem. 
\end{proof}
Putting the above lemmas together yields a proof of Theorem \ref{thm:Rice-Shapiro-theorem-for}. 
\begin{proof}
[Proof of Theorem \ref{thm:Rice-Shapiro-theorem-for}]Let $P$ be
a property of marked groups semi-decidable from recursive presentations.
Consider the set of finitely presented marked groups that belong to
$P$. This set is computably enumerable, since it is possible to enumerate
all finite presentations, and to select from these those that satisfy
$P$. 

All that we have to justify is that a group has $P$ if and only if
it is a quotient of one of these finitely presented marked groups. 

Suppose that a marked group $(H,S')$ is a strict marked quotient
of a group $(G,S)$ that has $P$. Consider the restriction of $P$
to $\{(H,S'),(G,S)\}$. It remains a $\nu_{RP}$-semi-decidable property.
It contains $(G,S)$. By Lemma \ref{lem:Two groups rec pres}, it
cannot be equal to $\{(G,S)\}$. Thus it must be all of $\{(H,S'),(G,S)\}$,
and $(H,S')$ also has $P$. 

Suppose now that $(G,S)$ is belongs to $P$. If $(G,S)$ is finitely
presentable, there is nothing to do. Otherwise, let $\mathcal{B}$
be the set of finitely presented marked groups that have $(G,S)$
as a quotient. By Lemma \ref{lem:Second-lemma:-sequence Rice Shapiro},
$\{(G,S)\}$ is not a $\nu_{RP}$-semi-decidable subset of $\{(G,S)\}\cup\mathcal{B}$,
and thus some group in $\mathcal{B}$ must belong to $P$. And thus
$(G,S)$ is indeed a quotient of a finitely presentable group with
$P$. 
\end{proof}

\section{\label{part: Co-rec presentations}Rice theorem for co-recursive
presentations}

One might expect to have a Rice-Shapiro theorem for groups given by
co-recursive presentations that mirrors exactly the Rice-Shapiro theorem
for recursive presentations, based on the exact same arguments, simply
using the lattice of marked groups with reverse order, and the corresponding
Scott topology. However, the proof of Theorem \ref{thm:Rice-Shapiro-theorem-for}
in fact relied crucially on the set of finitely presented groups,
which is a $\nu_{RP}$-computable sequence (all finite presentations
can be enumerated, and a finite presentation can be seen as a recursive
presentation, i.e. $\nu_{FP}\le\nu_{RP}$) of the compact elements
of $(\mathcal{G}_{k},\rightarrow)$.

Note however that there is no notion of ``finite co-presentation''.
For instance, there is no least marked group defined on a single generator
$a$ with $a^{2}\ne1$. Indeed, the marked groups $\langle a\,\vert\,a^{3}\rangle$
and $\langle a\,\vert\,a^{5}\rangle$ both satisfy that $a^{2}\ne1$,
and yet their $\inf$ in the lattice $(\mathcal{G}_{1},\rightarrow)$
is the trivial group. The closest notion is that of finite discriminating
family, which leads to the notion of \emph{absolute presentation}
of Neumann \cite{Neumann1973}, see \cite{Cornulier2007} for a modern
account. It was asked by Mann in \cite{Mann1982} whether it is possible
to algorithmically list all absolute presentations of groups, this
problem seems to still be open today. 

The following proposition also illustrates the fact that the meet
and join operations of $(\mathcal{G}_{k},\rightarrow)$ are not symmetrical
from the point of view of computability. 
\begin{prop}
\label{prop: Compute Meet and Join from presentations}The meet and
join operations of $(\mathcal{G}_{k},\rightarrow)$ are computable
for recursively presented groups.

The meet is computable for finitely presented groups, but the join
is not, since the join of two finitely presented groups does not have
to be finitely presented. Even: the join of two free groups does not
have to be finitely presented. 

The join operation is computable for co-recursively presented groups,
but the meet is not, since the meet of two co-recursively presented
groups does not have to be co-recursively presented. 
\end{prop}

\begin{proof}
\textbf{Meet operation.}

Consider two recursively presented groups given by $\pi_{1}=\langle S\vert R_{1}\rangle$
and $\pi_{2}=\langle S\vert R_{2}\rangle$. Their join is given by
$\pi_{1}=\langle S\vert R_{1},R_{2}\rangle$, which is indeed recursively
presented. The same holds modulo changing recursive for finite in
the above. 

Consider now a finitely presented group $\langle S\,\vert\,R\rangle$
with unsolvable word problem. Let $R=(r_{1},...,r_{n})$ and consider
the one-relator groups $\langle S\,\vert\,r_{i}\rangle$, for $1\le i\le n$.
These have solvable word problem by a theorem of Magnus. Consider
the meets $\langle S\,\vert\,r_{1}\rangle\wedge\langle S\,\vert\,r_{2}\rangle\wedge...\wedge\langle S\,\vert\,r_{i}\rangle$,
for each $i\le n$. For $i=n$, it has unsolvable word problem. The
least $i$ for which this is the case gives the meet of two groups
with solvable word problem that has unsolvable word problem. 

\textbf{Join operation.}

Consider two $k$-marked groups $(G,S)$ and $(H,S')$, and identify
$S$ and $S'$ via the canonical bijection. Let $\mathbb{F}_{S}$
be the free group over $S$, and express $(G,S)$ and $(H,S')$ as
marked quotients of the free group: $(G,S)=(\mathbb{F}_{S}/N_{1},S)$
and $(H,S)=(\mathbb{F}_{S}/N_{2},S)$. The meet $(G,S)\vee(H,S')$
is given by $(\mathbb{F}_{S}/(N_{1}\cap N_{2}),S)$. 

If $(G,S)$ and $(H,S')$ are recursively presented, then $N_{1}$
and $N_{2}$ are c.e. sets, and so is their intersection, thus $(G,S)\vee(H,S')$
is also recursively presented. And this is uniform. 

If $(G,S)$ and $(H,S')$ are co-recursively presented, then the complements
of $N_{1}$ and $N_{2}$ are c.e. sets, and so is their union, thus
$(G,S)\vee(H,S')$ is also co-recursively presented. And this is uniform. 

Finally, the fact that $(G,S)$ and $(H,S')$ can be finitely presented
without $(G,S)\vee(H,S')$ being so follows from well known results
on the Mihailova subgroup. Let $\langle s_{1},...s_{k}\,\vert\,r_{1},...,r_{n}\rangle$
be a finitely presentable group. Consider two $k+n$-markings of the
free group $\mathbb{F}_{S}$, namely 
\[
(\mathbb{F}_{S},(s_{1},...,s_{k},r_{1},...,r_{n})),
\]
\[
(\mathbb{F}_{S},(s_{1},...,s_{k},1,...,1)).
\]
The meet of these two markings is precisely the Mihailova subgroup
\cite{Mihailova1968}, which by a theorem of Grunewald \cite{Grunewald1978}
is not finitely presentable as soon as $\langle s_{1},...s_{k}\,\vert\,r_{1},...,r_{n}\rangle$
defines an infinite group.
\end{proof}
All this makes it so that we leave as an open problem the Rice-Shapiro
Theorem for co-recursive presentation. We still establish the Rice
Theorem. 
\begin{prop}
\label{prop: co-rec semi-D =00003D>upper stable  }If $P$ is a marked
group property that is semi-decidable from co-recursive presentations,
than any group that has a marked quotient with $P$ also has $P$. 
\end{prop}

This proposition follows directly from Lemma \ref{lem: co-rec pres pair of groups }.
As an immediate corollary we get:
\begin{cor}
[Rice theorem for co-recursive presentations]There is no non-trivial
decidable group property for groups given by co-recursive presentations. 
\end{cor}

\begin{proof}
Let $P$ be a $\nu_{co-RP}$-decidable property. Either $P$ or its
complement contains the trivial group. Then, by Proposition \ref{prop: co-rec semi-D =00003D>upper stable  },
the one that contains it is in fact the set of all groups. 
\end{proof}
The following lemma is essentially identical to Lemma \ref{lem:Two groups rec pres}.
\begin{lem}
\label{lem: co-rec pres pair of groups }Suppose that $(G,S)$ and
$(H,S')$ are two co-recursively presented marked groups, and that
$(H,S')$ is a strict marked quotient of $(G,S)$. Then no algorithm
that takes as input co-recursive presentations of either $(G,S)$
or $(H,S')$ can stop exactly on the presentations that define $(H,S')$. 
\end{lem}

\begin{proof}
Define a co-recursive presentation associated to a run of $\varphi_{n}(n)$:
enumerate the co-word problem of $(H,S')$ while this run last, if
it stops, start enumerating the co-word problem of $(G,S)$ instead.
\end{proof}
The following result is the equivalent of Lemma \ref{lem:Second-lemma:-sequence Rice Shapiro},
replacing recursive presentations by co-recursive presentations. Its
proof is identical to that of Lemma \ref{lem:Second-lemma:-sequence Rice Shapiro},
we thus omit it. The difference between this lemma and Lemma \ref{lem:Second-lemma:-sequence Rice Shapiro}
is that, here, we suppose that is given a computable sequence of co-recursively
presented groups whose supremum is $(G,S)$. In Lemma \ref{lem:Second-lemma:-sequence Rice Shapiro},
such a sequence (with reversed arrows) was constructed from $(G,S)$
by taking truncated presentations. A similar construction -``truncated''
co-recursive presentations- is not available here. 
\begin{lem}
Suppose that $(G,S)$ is a co-recursively presented marked group.
Let $(G_{n},S_{n})$ be a $\nu_{co-RP}$-computable sequence of marked
quotients of $(G,S)$, with for each $n$ $(G_{n+1},S_{n+1})\rightarrow(G_{n},S_{n})$,
and whose supremum is $(G,S)$. 

Then no algorithm that takes as input co-recursive presentations of
marked groups in $\{(G,S)\}\cup\{(G_{n},S_{n}),n\in\mathbb{N}\}$
can stop exactly on the presentations that define $(G,S)$. 
\end{lem}

\begin{problem}
[Rice-Shapiro for co-recursive presentations]Is it true that if
$P$ is a property of marked groups that is semi-decidable from co-recursive
presentations, then there exists a computable sequence $(\mathcal{A}_{n})_{n\in\mathbb{N}}$
of finite sets of relations, such that a marked group $(G,S)$ satisfies
$P$ if and only if there exists $n$ so that no relation of $\mathcal{A}_{n}$
holds in $(G,S)$? 
\end{problem}

\section{decision problem from word problem algorithms}

By taking the conjunction $\nu_{RP}\wedge\nu_{co-RP}$, we obtain
the numbering associated to word problem algorithms, which we call
$\nu_{WP}$. 

In \cite{Rauzy2021b}, we investigate in details decision problems
for groups given by word problem algorithms. The topology that is
relevant to study such decision problems is the topology of \emph{the
space of marked groups}, which can be defined thanks to a metric defined
as follows: two $k$-marked groups are at distance less or equal to
$2^{-n}$ if and only if they satisfy exactly the same relations of
length at most $n$. 

We want to note here that the relationship between the topology of
the space of marked groups and decision problems for groups given
by word problem algorithms is looser than the one that relates problems
about recursive presentations and the Scott topology. Indeed, it is
false that a property semi-decidable for $\nu_{WP}$ is open in the
space of marked groups (by a result of Friedberg about the computable
reals, see for instance \cite{Hoyrup2016}). It is still possible
that the decidable properties of $\nu_{WP}$ are automatically clopen
in the space of marked groups, this is an open instance of what is
known as the \emph{continuity problem} \cite{SPREEN2016}. However,
we can show (by \cite{Moschovakis1964}) that a computable sequence
which is dense in a certain open set of the space of marked groups
must also meet every $\nu_{WP}$-semi-decidable set. This indicates
that the topology of the space of marked groups is close to the topology
of $\nu_{WP}$-semi-decidable sets. 

\section{\label{sec:Marked-quotient-algorithms-FP}Marked quotient algorithms
and finite presentations}

\subsection{Algorithmic characterization of finitely presented groups }

Recall from the introduction that a marked group $(G,S)$ has a \emph{marked
quotient algorithm} if there is a procedure that stops exactly on
recursive presentations of its marked quotients, i.e. if the set $\{(H,S'),\,(G,S)\rightarrow(H,S')\}$
is $\nu_{RP}$-semi-decidable. 
\begin{lem}
\label{lem:Any-finitely-presented-has-marked-quotients}Any finitely
presented group admits a marked quotient algorithm. And this is uniform. 
\end{lem}

\begin{proof}
To check whether the relation $(G,S)\rightarrow(H,S')$ holds, it
suffices to check whether the finitely many relations of $(G,S)$
hold in $(H,S')$, this can be semi-decided thanks to a recursive
presentation. 
\end{proof}
Note that any non recursively presentable simple group has a marked
quotient algorithm. When restricting our attention to recursively
presented groups, having a marked quotient algorithm will characterize
finitely presentable groups. Recall that $\nu_{MQA}$ is the numbering
of finitely presented groups associated to marked quotient algorithms. 
\begin{thm}
\label{thm:AlgorithmicCharacterizationFPgroups}A group is finitely
presented if and only if it admits both a recursive presentation and
a marked quotient algorithm. And this statement is uniform: it in
fact provides an equivalence of numberings:
\[
\nu_{FP}\equiv\nu_{RP}\wedge\nu_{MQA}.
\]
\end{thm}

\begin{proof}
Lemma \ref{lem:Any-finitely-presented-has-marked-quotients} provides
the obvious direction. 

Suppose now that we have a recursive presentation $\langle S\,\vert\,r_{1},r_{2},...\rangle$
and a marked quotient algorithm for a marked group $(G,S)$. Let $(G_{n},S)$
be the marked group given by the truncated presentation $\langle S\vert r_{0},r_{1},r_{2},...,r_{n}\rangle$.

By Lemma \ref{lem:Second-lemma:-sequence Rice Shapiro}, no algorithm
that accepts recursive presentations of groups in $\{(G,S)\}\cup\{(G_{n},S),\,n\in\mathbb{N}\}$
can stop exactly on those that define $(G,S)$. 

But the marked quotient algorithm of $(G,S)$ does accept every presentation
that defines $(G,S)$, since $(G,S)$ is a marked quotient of itself.
Thus there must be another group in $\{(G,S)\}\cup\{(G_{n},S),\,n\in\mathbb{N}\}$
which is a marked quotient of $(G,S)$. But for each $n$ we have
a morphism of marked groups $(G_{n},S)\rightarrow(G,S)$. Thus if
a morphism exists in the other direction, it must be that $(G,S)$
is isomorphic to some $(G_{n},S)$, and thus it is finitely presented.

We now have to justify that the above proof is effective, in that
some $n$ with $(G,S)=(G_{n},S)$ can be computed from the recursive
presentation of $(G,S)$ and its marked quotient algorithm.

Going back to the proof of Lemma \ref{lem:Second-lemma:-sequence Rice Shapiro},
we see that was constructed a computable sequence of recursive presentations
$(E_{n})_{n\in\mathbb{N}}$, such that if $\varphi_{n}(n)\uparrow$,
$E_{n}$ is the infinite presentation of $(G,S)$, and if $\varphi_{n}(n)\downarrow$,
$E_{n}$ defines $(G_{n},S)$. 

The sets $\{n\in\mathbb{N},\,\varphi_{n}(n)\downarrow\}$ and $\{n\in\mathbb{N},(G,S)\rightarrow E_{n}\}$
are two recursively enumerable sets, by the above non-effective argument
they intersect, thus it is possible to algorithmically produce an
element in their intersection. This will precisely provide a finite
presentation of $(G,S)$.
\end{proof}

\subsection{Decision problems from marked quotient algorithms }

In Section \ref{part:Rice-Shapiro}, we established the Rice-Shapiro
theorem for recursive presentations, which shows that very little
can be said of groups given by recursive presentations. We have now
decomposed the numbering associated to finite presentations into a
disjunction, $\nu_{FP}\equiv\nu_{RP}\wedge\nu_{MQA}$. We will now
show that this disjunction is non-trivial, we thus must show that
$\nu_{FP}\not\equiv\nu_{MQA}$ (since it follows from Section \ref{part:Rice-Shapiro}
that $\nu_{FP}\not\equiv\nu_{RP}$). We will thus show that very little
can be said about a group given by a marked quotient algorithm. Note
that in the following, we only consider finitely presented groups,
but these are given by descriptions that are much weaker than finite
presentations. 
\begin{lem}
\label{lem:Two groups MQAlgo}Suppose that $(G,S)$ and $(H,S')$
are two finitely presented marked groups, and that $(H,S')$ is a
strict marked quotient of $(G,S)$. Then no algorithm that takes as
input marked quotient algorithms of either $(G,S)$ or $(H,S')$ can
stop exactly on those that define $(H,S')$. 
\end{lem}

In other words, $\{(G,S)\}$ is not a $\nu_{MQA}$-semi-decidable
subset of $\{(G,S),(H,S')\}$. 
\begin{proof}
The proof is similar to that of Lemma \ref{lem:Two groups rec pres},
we only sketch it. Associated to a run of $\varphi_{n}(n)$, construct
an algorithm that recognizes the quotients of $(H,S')$ while this
run last, and that also accepts the quotients of $(G,S)$ if this
run stops.
\end{proof}
Call a property $P$ of marked groups \emph{upward closed} when a
marked group that has a quotient in $P$ also has $P$.
\begin{cor}
[Rice-Shapiro for marked quotient algorithms] The $\nu_{MQA}$-semi-decidable
properties are exactly the $\nu_{FP}$-semi-decidable properties that
are upward closed in $(\mathcal{G},\rightarrow)$. 
\end{cor}

\begin{proof}
Let $P$ be a property of finitely presented marked groups that is
$\nu_{MQA}$-semi-decidable. That it is semi-decidable from finite
presentations is immediate: finite presentations provide strictly
more information than marked quotient algorithms. That it is upward
closed follows from Lemma \ref{lem:Two groups MQAlgo}. 

Suppose now that $P$ is semi-decidable from finite presentations.
Then, the property ``having a marked quotient in $P$'' is semi-decidable
from marked quotient algorithms, since to check if a marked group
$(G,S)$ has a marked quotient in $P$, it suffices to enumerate all
finite presentations of groups in $P$, and to apply to all of these
in parallel the marked quotient algorithm of $(G,S)$. When $P$ is
upward closed, ``having a marked quotient in $P$'' is simply $P$.
\end{proof}

\section{\label{sec:Marked-quotient-algorithms-Relative}Marked quotient algorithms
relative to a class of groups}

In this section, we introduce relative marked quotient algorithms,
and explain that the cases where these are most interesting is when
they do not rely on notions of finite presentations modulo a certain
class of groups. 

\subsection{Marked quotients algorithms and semi-decidable equality}
\begin{defn}
\label{def:RELATIVEQALGO}Let $\mathcal{C}$ be a class of marked
groups. We say that a marked group $(G,S)$ has \emph{a marked quotient
algorithm relative to $\mathcal{C}$} if there is an algorithm that
takes as input recursive presentations for groups in $\mathcal{C}$
and stops exactly on those that define marked quotients of $(G,S)$. 
\end{defn}

See Proposition \ref{prop:The-Lamplighter-group-CFQ} for an example
of a group that has a relative marked quotient algorithm for finite
groups without being finitely presented. 

One easily defines a numbering $\nu_{MQA}^{\mathcal{C}}$ associated
to\emph{ $\mathcal{C}$}-marked quotient algorithms, modifying the
definition of $\nu_{MQA}$ accordingly (see Definition \ref{def:Marqued quotient algorithm}).
(In cases where a group is not uniquely determined by its marked quotients
in $\mathcal{C}$, when it is not residually $\mathcal{C}$, the marked
quotient algorithm in $\mathcal{C}$ may not determine a marked group
uniquely. In this case, $\nu_{MQA}^{\mathcal{C}}$ is a multi-numbering
of marked groups, i.e. a partial multi-function $\nu_{MQA}^{\mathcal{C}}:\subseteq\mathbb{N}\rightrightarrows\mathcal{G}$.
We will however not need that notion here.)

Note that if we suppose that $\mathcal{C}$ is a class of groups,
i.e. is closed under isomorphism, then ``having a marked quotient
relative to $\mathcal{C}$'' becomes a group property.
\begin{lem}
\label{lem:Change of marking }Let $\mathcal{C}$ be a class of abstract
groups. If a group $G$ admits a marked quotient algorithm relative
to $\mathcal{C}$ with respect to some marking, then it also admits
one with respect to any marking.
\end{lem}

\begin{proof}
Let $S$ and $T$ be two finite generating sets of a group $G$. We
suppose that we have access only to the marked quotient algorithm
for $G$ with respect to $S$. Fix for each $s$ in $S$ an expression
$s=t_{1}^{\alpha_{1}}...t_{k}^{\alpha_{k}}$, with $\alpha_{i}\in\left\{ -1,1\right\} $
and $t_{i}\in T$, that gives $s$ as a product of elements of $T$
and of their inverses, and for each $t$ in $T$ an expression $t=s_{1}^{\beta_{1}}s_{2}^{\beta_{2}}...s_{k}^{\beta_{k}}$
that describes $t$ in terms of the generators of $S$ and their inverses. 

Consider a marked group ($H,T')$ given by a recursive presentation,
with $T'$ of same cardinality as $T$. We want to determine whether
$(H,T')$ is a marked quotient of $(G,T)$, thus whether the bijection
$f:T\rightarrow T'$ can be extended to a group morphism $\tilde{f}$. 

Let $S'$ be the family defined, in $H$, by the same formulas as
$S$ is in $G$, i.e. for each $s=t_{1}^{\alpha_{1}}...t_{k}^{\alpha_{k}}$
in $S$, we define $s'=f(t_{1})^{\alpha_{1}}...f(t_{k})^{\alpha_{k}}$.
Notice that if $f$ does define a morphism $\tilde{f}$ of $G$ onto
$H$, $S'$ should be the image of the family $S$ by $\tilde{f}$,
and thus it should be a generating family of $H$. We can therefore,
using the recursive presentation of $(H,T')$, look for an expression
of the elements of $T'$ in terms of the elements of $S'$ in $H$. 

If such an expression does not exist, our procedure will not stop,
but then $(H,T')$ is not a quotient of $(G,T)$, thus this result
is coherent. 

Otherwise we can use the formulas just found to obtain a recursive
presentation of $H$ with respect to $S'$. From it, we can ask whether
the natural bijection $S\rightarrow S'$ defines a group morphism,
thanks to the marked quotient algorithm of $G$ on $S$. If this procedure
does not end, $(H,T')$ is not a marked quotient of $(G,T)$. 

If it terminates, the bijection $S\rightarrow S'$ does define a group
morphism $\psi:G\rightarrow H$. But we still have to check that $\psi$
maps $T$ to $T'$. This is done using the expressions of the form
$t=s_{1}^{\beta_{1}}s_{2}^{\beta_{2}}...s_{k}^{\beta_{k}}$, that
define, in $G$, the elements of $T$ in term of those of $S$. 

Compute the formal images of these elements inside of $H$. This yields
for each $t$ in $T$ an expression of $\psi(t)$ in terms of the
elements of $S'$. 

It then follows that $f$ extend to a group homomorphism if and only
if for each $t$ of $T$, $\psi(t)=f(t)$ in $H$. This can be semi-decided
thanks to the recursive presentation of $H$. 
\end{proof}
The following is an easy consequence of the fact that the quotient
relation is an order for isomorphism classes of marked groups: 
\begin{fact}
Isomorphism of marked groups is semi-decidable for groups in $\mathcal{C}$
that are given by $\nu_{RP}\wedge\nu_{MQA}^{\mathcal{C}}$-names.
\end{fact}

\begin{proof}
Given $\nu_{RP}\wedge\nu_{MQA}^{\mathcal{C}}$-names for two marked
groups, it suffices to ask whether the first is a marked quotient
of the second, and whether the second is a marked quotient of the
first one.
\end{proof}
\begin{rem}
The above fact can be extended to semi-decidability of the actual
isomorphism relation when $\mathcal{C}$ is closed under group isomorphism,
by noticing that Lemma \ref{lem:Change of marking } is effective. 

Having a semi-decidable isomorphism problem seems to be a very good
criterion to ensure that a certain numbering $\nu$ is a good candidate
for Problem \ref{prob: Descriptions that are strong enough}. In particular,
it guarantees that the isomorphism problem will always be solvable
on \emph{finite} sets of groups. It also guarantees that precise classification
theorems (such as that of finitely generated abelian groups) immediately
translate into an effective result. (And thus while the proof of the
classification of finitely generated abelian groups is effective,
the statement of this theorem itself is sufficient to argue that any
numbering $\nu$ of finitely generated abelian groups, which provides
less information than the finite presentation numbering, i.e. for
which $\nu\ge\nu_{FP}$, and which has a semi-decidable isomorphism
problem, will permit to actually solve the isomorphism problem for
finitely generated abelian groups.)
\end{rem}

\subsection{Finite presentation modulo a class of groups }

Theorem \ref{thm:AlgorithmicCharacterizationFPgroups} states that
being finitely presentable is equivalent to being recursively presentable
and having a marked quotient algorithm. 

A similar phenomenon sometimes occur for relative marked quotient
algorithms, when we replace finite presentations by other appropriate
notions of ``finite presentation relative to a class of groups''. 

In fact, the introduction of relative quotient algorithms will be
most interesting when there is no underlying notion of relative finite
presentation. 

A marked group $(G,S)$ is called residually $\mathcal{C}$ when any
non-trivial element of it has a non-trivial image in a marked quotient
that belongs to $\mathcal{C}$. Every marked group has a greatest
(for the order of the lattice of marked groups) marked quotient which
is residually $\mathcal{C}$, called its residually $\mathcal{C}$
image. 

For $(G,S)$ a marked group and $\mathcal{C}$ a set of marked groups,
let $\mathcal{E}((G,S),\mathcal{C})$ be the set of morphisms that
exist between $(G,S)$ and groups in $\mathcal{C}$.
\begin{defn}
Let $\mathcal{C}$ be a class of marked groups, and $(G,S)$ a marked
group. The \emph{residually $\mathcal{C}$ image} of $(G,S)$ is the
group 
\[
G\big/\bigcap_{\phi\in\mathcal{E}((G,S),\mathcal{C})}\ker(\phi).
\]
\end{defn}

Thanks to the above, we can define finite presentations modulo $\mathcal{C}$. 
\begin{defn}
A residually $\mathcal{C}$ marked group $(G,S)$ has a \emph{finite
presentation as a residually $\mathcal{C}$ group }if there exists
a finitely presented group whose residually $\mathcal{C}$ image is
$(G,S)$. 
\end{defn}

The following provides the obvious case for when a marked group has
a marked quotient algorithm relative to a class $\mathcal{C}$. 
\begin{prop}
Suppose that the residually $\mathcal{C}$ image of $(G,S)$ is finitely
presentable as a residually $\mathcal{C}$ group. Then $(G,S)$ has
a marked quotient relative to $\mathcal{C}$.
\end{prop}

\begin{proof}
Just as Lemma \ref{lem:Any-finitely-presented-has-marked-quotients}. 
\end{proof}

\begin{example}
By equational noetherianity of free groups, for every marked group
$(G,S)$, there is a finitely presented group $(H,S')$ such that
$(G,S)$ and $(H,S')$ have exactly the same marked free quotients.
See for instance \cite[Theorem 2.7]{Houcine2007}. Thus every group
has a marked quotient algorithm for free groups. Similarly, every
group has a marked quotient algorithm relative to abelian groups.
\end{example}

We can also show that with respect to group varieties, all marked
quotient algorithms in fact rely on finite presentations. If $\mathcal{V}$
is a group variety, a residually $\mathcal{V}$ group is simply a
group in $\mathcal{V}$, because group varieties are closed under
unrestricted direct products and subgroups. And a finite presentation
as a residually $\mathcal{V}$ group is simply a finite presentation
which omits the infinitely many relations that are given by the laws
of the variety. We call these \emph{finite presentations modulo $\mathcal{V}$}
rather than finite presentation of residually $\mathcal{V}$ groups.
\begin{prop}
\label{prop: Quotient algo for group variety }Let $\mathcal{V}$
be a group variety defined by a c.e. set of laws. A recursively presented
marked group $(H,S)$ of \textup{$\mathcal{V}$} admits a \textup{$\mathcal{V}$}-marked
quotient algorithm if and only if is finitely presented modulo $\mathcal{V}$
. 

What's more, this statement is uniform: it states that the numbering
associated to finite presentations modulo $\mathcal{V}$ is equivalent
to the conjunction $\nu_{RP}\wedge\nu_{MQA}^{\mathcal{V}}$.
\end{prop}

\begin{proof}
The proof is identical to that of Theorem \ref{thm:AlgorithmicCharacterizationFPgroups},
except that one needs to add the laws that defines $\mathcal{V}$
when building the presentations $E_{n}$.
\end{proof}

\subsection{A relative marked quotient algorithm that does not rely on a finite
presentation }

We end this section by giving an example of a marked quotient algorithm
that does not rely on a finite presentation. The fact that the lamplighter
group has a marked finite quotient algorithm follows from results
of Hartung \cite{HARTUNG2011}, but we give a short proof as an illustration.
\begin{prop}
\label{prop:The-Lamplighter-group-CFQ}The Lamplighter group $\mathbb{Z}/2\mathbb{Z}\wr\mathbb{Z}$
has a marked quotient algorithm with respect to the set of finite
groups. 
\end{prop}

\begin{proof}
Indeed, it admits the following presentation: 
\[
\langle a,\varepsilon\,\vert\,\varepsilon^{2},\,\left[\varepsilon,a^{-n}\varepsilon a^{n}\right],n\in\mathbb{Z}\rangle
\]
To see whether a finite group $F$ generated by two elements $a_{1}$
and $\varepsilon_{1}$ is a quotient of it, find a multiple $N$ of
the order of $a_{1}$ using the recursive presentation of $F$. Then,
notice that $(F,(a_{1},\epsilon_{1}))$ is a quotient of $(L,(a,\epsilon))$
if and only if it is a quotient of the group obtained from $(L,(a,\epsilon))$
by adding the relation $a^{N}$. But the quotient $L/\langle\langle a^{N}\rangle\rangle$
is in fact the finite wreath product $\mathbb{Z}/2\mathbb{Z}\wr\mathbb{Z}/N\mathbb{Z}$,
which admits the finite presentation: 
\[
\langle a,\varepsilon\,\vert\,\varepsilon^{2},\,a^{N},\,\left[\varepsilon,a^{-n}\varepsilon a^{n}\right],0\leq n\leq N\rangle
\]
Thus the problem is reduced to checking finitely many relations in
$(F,(a_{1},\epsilon_{1}))$, this is semi-decidable.
\end{proof}
\begin{prop}
The lamplighter group is not finitely presented as a residually finite
group.
\end{prop}

\begin{proof}
Because the lamplighter group is residually finite, we just have to
prove that no finitely presented group has $\mathbb{Z}/2\mathbb{Z}\wr\mathbb{Z}$
as its residually finite image. 

Suppose that we have a marked group $(G,(a_{1},\epsilon_{1}))$, given
by a presentation $\langle a_{1},\varepsilon_{1}\vert r_{1},...,r_{p}\rangle$,
and a morphism $\phi:(G,(a_{1},\epsilon_{1}))\rightarrow(\mathbb{Z}/2\mathbb{Z}\wr\mathbb{Z},(a,\epsilon))$
which satisfies that any morphism $h$ from $(G,(a_{1},\epsilon_{1}))$
to a finite group $(F,(a_{2},\epsilon_{2}))$ factors through $\phi$.

Since $(\mathbb{Z}/2\mathbb{Z}\wr\mathbb{Z},(a,\epsilon))$ is a quotient
of $(G,(a_{1},\epsilon_{1}))$, it must satisfy the finitely many
relations of $(G,(a_{1},\epsilon_{1}))$. These relations must in
turn be consequences of a finite number of the relations of $(\mathbb{Z}/2\mathbb{Z}\wr\mathbb{Z},(a,\epsilon))$.
In particular, there must be a natural number $N$ such that the first
$N$ relations of $(\mathbb{Z}/2\mathbb{Z}\wr\mathbb{Z},(a,\epsilon))$
imply those of $G$. Consider $(H,(a,\varepsilon))$ given by the
presentation 
\[
\langle a,\varepsilon\,\vert\,\varepsilon^{2},\,\left[\varepsilon,a^{-n}\varepsilon a^{n}\right],0\le n\le N\rangle.
\]

Then $(G,(a_{1},\epsilon_{1}))$, $(H,(a,\varepsilon))$ and $(\mathbb{Z}/2\mathbb{Z}\wr\mathbb{Z},(a,\epsilon))$
must have exactly the same marked finite quotients. 

\begin{center}
	\begin{tikzcd}
	G \arrow[r, "\phi_{1}"]\arrow[rrd, "h_{0}"]& H \arrow[r, "\phi_{2}"]\arrow[rd, "h_{1}"]& L \arrow[ d, "h_{2}"]\\ && F
\end{tikzcd} 
\end{center}

To end the proof, we find a contradiction, by finding a finite group
which satisfies the relations of $(H,(a,\varepsilon))$, but not those
of $(\mathbb{Z}/2\mathbb{Z}\wr\mathbb{Z},(a,\epsilon))$. We define
a subgroup of the group $\mathfrak{S}_{5N}$ of permutations on $\{1,...,5N\}$.
Consider the element $\sigma_{0}$ of $\mathfrak{S}_{5N}$, defined
by the following formula:
\[
\sigma_{0}(i)=\begin{cases}
i+2 & i\le5N-2\\
i+2-5N & i\ge5N-1
\end{cases}.
\]
Let $\sigma_{1}$ be the product of the transpositions $(1,2)$ and
$(2N+4,2N+5)$. It is then easy to see that the following relations
hold between $\sigma_{0}$ and $\sigma_{1}$:
\[
\sigma_{1}^{2}=id,
\]
\[
\left[\sigma_{1},\sigma_{0}^{-n}\sigma_{1}\sigma_{0}^{n}\right]=id,\,1\le n\le N,
\]
\[
\left[\sigma_{1},\sigma_{0}^{-N-1}\sigma_{1}\sigma_{0}^{N+1}\right]\ne id.
\]
The subgroup of $\mathfrak{S}_{5N}$ generated by $\sigma_{0}$ and
$\sigma_{1}$ is thus a marked quotient of $(H,(a,\varepsilon))$,
but not of $(\mathbb{Z}/2\mathbb{Z}\wr\mathbb{Z},(a,\epsilon))$.
This contradicts the supposition that $(\mathbb{Z}/2\mathbb{Z}\wr\mathbb{Z},(a,\epsilon))$
is finitely presented as a residually finite group. 
\end{proof}
The following problem asks for an instance of a set of marked groups
with a semi-decidable equality that comes from a genuinely computational
generalization of the notion of finite presentation. 
\begin{problem}
Find a class $\mathcal{C}$ of recursively presented groups such that
all groups in $\mathcal{C}$ admit marked quotient algorithms relative
to $\mathcal{C}$ (so that equality is $\nu_{RP}\wedge\nu_{MQA}^{\mathcal{C}}$-semi-decidable)
but such that not all groups in $\mathcal{C}$ are finitely presented
as residually $\mathcal{C}$ groups.
\end{problem}

\section{\label{sec:The-Adian-Rabin-theorem}On the Adian-Rabin theorem }

In this section, we remark how the Adian-Rabin provides incomplete
answers to certain natural problems. We first recall its proof. 
\begin{defn}
A \emph{Markov property} is a group property $P$ which admits a \emph{positive
witness}, any finitely presented group with $P$, and a \emph{negative
witness}, which is any finitely presented group that does not embed
in any group with $P$. 
\end{defn}

\begin{thm}
[Adian-Rabin]A Markov property is $\Sigma_{1}^{0}$-hard, and thus
not co-semi-decidable, for groups given by finite presentations. 
\end{thm}

We sketch the proof given in \cite{MillerIII1992}, omitting the main
technical lemma, which is the actual content of the proof, but our
purpose is to render explicit what is proved exactly. 
\begin{proof}
Fix $G_{+}$ and $G_{-}$ the positive and negative witnesses of a
Markov property. 

Firstly, given a finitely presented marked group $(H,S)$ with unsolvable
word problem, a family of finite presentations $\pi_{w}$, $w\in(S\cup S^{-1})^{*}$,
is constructed, such that: $\pi_{w}$ can be effectively constructed
from $w$, and $\pi_{w}$ defines the trivial group if and only if
$w=1$ in $H$. Whatever $w$, $\pi_{w}$ is defined on two generators
$u$ and $v$. Finally, when $w\ne1$ in $(H,S)$, $\pi_{w}$ contains
a copy of $G_{-}$. 

Consider a presentation for the free product of $G_{+}$ and of the
group defined by $\pi_{w}$. 

This defines a group that contains $G_{-}$ when $w\ne1$, and on
the contrary this defines exactly $G_{+}$ when $w=1$ in $(H,S)$. 
\end{proof}

\subsection{\label{subsec: marked group recognition adian rabin }From the point
of view of marked groups }

Consider the problem of recognizing a single marked group from a finite
presentation. 

This problem is always semi-decidable. 

This problem is decidable for free groups, when the marking is via
a basis: a finite presentation defines a free group marked by a basis
if and only if it has no relations modulo free reductions. It is also
decidable for cyclic groups marked by families with one generators. 

What we can easily extract from the proof of the Adian-Rabin is the
following: 
\begin{prop}
\label{prop:UNRECON-1}Consider a finitely presented marked group
$(G,(s_{1},...,s_{k}))$. The problem of deciding whether a finite
presentation defines the marked group $(G,(s_{1},...,s_{k},1_{G},1_{G}))$
is not co-semi-decidable. 
\end{prop}

\begin{proof}
This follows immediately from the proof of the Adian-Rabin Theorem
described above, by taking $G_{+}=G$. The two additional generators
corresponding to the identity appear in the free product with the
presentation $\pi_{w}$.
\end{proof}
In \cite{Baumslag2007}, Baumslag and Miller establish the following
result:
\begin{thm}
[\cite{Baumslag2007}, Theorem B]If n is chosen sufficiently large,
then there is no algorithm to determine of a finite presentation with
$n$ generating symbols whether or not the group presented is free
abelian of rank $n$.
\end{thm}

This is, to the best of our knowledge, the only result that provides
more information than the Adian-Rabin construction about recognition
of marked groups thanks to finite presentations. 
\begin{problem}
Bridge the ``two additional generators gap'' that appears in Proposition
\ref{prop:UNRECON-1}, and describe exactly the marked groups that
define decidable properties of groups given by finite presentations. 
\end{problem}

\subsection{Classification of properties in the arithmetical hierarchy}

The Adian-Rabin theorem does not permit to precisely classify properties
of finitely presented groups in the arithmetical hierarchy, except
at the very first levels. 

Most of the group properties whose precise location in the arithmetical
hierarchy is known are properties that are semi-decidable but not
co-semi-decidable. The Adian-Rabin Theorem is used to prove the ``not
co-semi-decidable'' part, and semi-decidability is proved directly.
This includes: being abelian, trivial, finite, nilpotent, virtually
nilpotent, hyperbolic, having Kazhdan's Property (T) \cite{Ozawa2014}.
Note that ``having a non-trivial finite quotient'' is semi-decidable
and not decidable, but that this does not follow from the Adian-Rabin
Theorem \cite{Bridson2015b}. 

However, when it comes to classifying properties that are higher in
the arithmetical hierarchy, almost nothing is known. In particular,
the group properties one naturally encounters are usually Markov or
co-Markov properties, \emph{but not both. }And thus the Adian-Rabin
falls short. 

``Having solvable word problem'' is known to be $\Sigma_{3}^{0}$-complete
\cite{Boone1966}. 

``Being torsion free'' is known to be $\Pi_{2}^{0}$-complete, as
was shown by Lempp \cite{Lempp1997}. The proof of this result amounts
to noticing that the Higman Embedding Theorem preserves torsion-freeness,
it is thus easy to adapt to other properties that known Higman embeddings
preserve, like ``having a non-trivial center'' which is $\Sigma_{2}^{0}$-complete
\cite{OuldHoucine2007}, but hard to adapt to other properties. (For
instance, establishing that being solvable is not semi-decidable is
probably easier than establishing a Higman Embedding Theorem for soluble
groups.) 

Note also that, as a corollary of \cite{Baumslag1983}, it is proved
in \cite{Downey2008} that deciding whether two given finitely presented
groups have the same integral homology sequence is $\Sigma_{1}^{1}$-complete
-thus not even in the arithmetical hierarchy.

We gave a proof in Proposition \ref{prop: RF not CE} of the fact
that ``being residually finite'' is not semi-decidable. Note that
this property is Markov, so it was already known to be not co-semi-decidable,
but it is not co-Markov. 

The following problem goes back to Miller \cite{MillerIII1992} and
Lempp \cite{Lempp1997}. 
\begin{problem}
Classify group properties higher up in the arithmetical hierarchy,
for groups given by finite presentations. Prove that the following
properties are not semi-decidable from finite presentations: being
solvable, or being $k$-solvable for a fixed $k\ge2$, being simple,
being torsion, being amenable. 
\end{problem}

Note that the Adian-Rabin Theorem provides a result similar to \emph{half}
of the Rice-Shapiro Theorem: it corresponds to Lemma \ref{lem:Two groups rec pres}
adapted to groups given by finite presentations. No equivalent of
the other half of the Rice-Shapiro Theorem, Lemma \ref{lem:Second-lemma:-sequence Rice Shapiro},
was ever proved for groups given by finite presentations. In other
words, we can make statements such as: ``this property is not semi-decidable,
because checking it would involve proving that a relation is not satisfied
in the given group, which is impossible'', but we cannot make statements
such as ``this property is not semi-decidable, because checking it
would involve proving simultaneously that infinitely many independent
relations are satisfied in the given group, which is impossible''.
Although presumably, at least for some of the properties quoted above,
this will be true. 

\subsection{Classifying properties for groups given by finite presentations and
word problem algorithms}

The proof given above of the Adian-Rabin Theorem in fact relies on
a single marking of the positive witness, the group $G_{+}$. We already
noted that this marking must be extended by twice the identity as
new generators for the construction to work. What we want to note
here is that the fact that this construction relies on a single marking
of $G_{+}$ makes it so that it provides no information on decision
problems for groups given by word problem algorithms and recursive
presentations. 

Indeed, we have the obvious fact:
\begin{prop}
The marked isomorphism problem is decidable for $\nu_{WP}\wedge\nu_{FP}$.
\end{prop}

\begin{proof}
To check whether two finite presentations $\langle S\,\vert\,R\rangle$
and $\langle S\,\vert\,R'\rangle$ define the same marked group, it
suffices to decide whether the relations of $R$ hold in $\langle S\,\vert\,R'\rangle$,
and if those of $R'$ hold in $\langle S\,\vert\,R\rangle$. This
can be decided if we can solve the word problem in these groups. 
\end{proof}
And thus any result that establishes that a certain property $P$
is not $\nu_{WP}\wedge\nu_{FP}$-semi-decidable will have to rely
on using infinitely many different markings of groups without $P$,
and infinitely many markings of groups with $P$. The Adian-Rabin
construction uses infinitely many markings of groups without $P$,
but only a single marking of a group with $P$.

For instance, in \cite[Theorem 26, Chapter IV]{MillerIII1972}, Miller
proved that the isomorphism problem is unsolvable for finitely presented
residually finite groups by building a finitely presented residually
finite group $G$ such that the set $[G]$ of all markings of $G$
is $\nu_{FP}$-semi-decidable but not $\nu_{FP}$-decidable \emph{amongst
residually finite groups}. Because the word problem is uniformly solvable
for finitely presented residually finite groups, the marked isomorphism
problem is solvable for these. Thus a necessary feature of Miller's
construction is that infinitely many markings of $G$ are used. Even
more: fixing any generating family $S$ of $G$, the average word
length with respect to $S$ of the different generating families of
$G$ that are used grows faster than every recursive functions.

We thus end this paragraph by the following problem: 
\begin{problem}
Establish a theorem that gives sufficient conditions on a group property
for it to be undecidable for groups given by finite presentations
and word problem algorithms. 
\end{problem}

\bibliographystyle{alpha}
\bibliography{AlgorithmicBiblio2}

\newcommand{\etalchar}[1]{$^{#1}$}
\begin{thebibliography}{BBMZ23}

\bibitem[Bar03]{Bartholdi2003}
Laurent Bartholdi.
\newblock Endomorphic presentations of branch groups.
\newblock {\em Journal of Algebra}, 268(2):419--443, oct 2003.

\bibitem[BBMZ23]{Belk2023}
James Belk, Collin Bleak, Francesco Matucci, and Matthew C.~B. Zaremsky.
\newblock Progress around the boone-higman conjecture.
\newblock {\em arxiv:2306.16356}, June 2023.

\bibitem[BDM83]{Baumslag1983}
G.~Baumslag, E.~Dyer, and C.F. Miller.
\newblock On the integral homology of finitely presented groups.
\newblock {\em Topology}, 22(1):27--46, 1983.

\bibitem[BEH08]{Bartholdi2008}
Laurent Bartholdi, Bettina Eick, and Ren{\'{e}} Hartung.
\newblock A nilpotent quotient algorithm for certain infinitely presented
  groups and its applications.
\newblock {\em International Journal of Algebra and Computation},
  18(08):1321--1344, 12 2008.

\bibitem[BKM07]{Bumagin2007}
Inna Bumagin, Olga Kharlampovich, and Alexei Miasnikov.
\newblock The isomorphism problem for finitely generated fully residually free
  groups.
\newblock {\em Journal of Pure and Applied Algebra}, 208(3):961--977, March
  2007.

\bibitem[BM37]{Banach1937}
Stefan Banach and Stanislaw Mazur.
\newblock Sur les fonctions calculables.
\newblock {\em Ann. Soc. Pol. de Math}, 16(223):402, 1937.

\bibitem[BMI07]{Baumslag2007}
Gilbert Baumslag and Charles~F. Miller~III.
\newblock The isomorphism problem for residually torsion-free nilpotent groups.
\newblock {\em Groups, Geometry, and Dynamics}, 1(1):1--20, March 2007.

\bibitem[BR66]{Boone1966}
William~W. Boone and Hartley~Jr. Rogers.
\newblock On a problem of {J.H.C.} {W}hitehead and a problem of {A}lonzo
  {C}hurch.
\newblock {\em Mathematica Scandinavica}, 19:185, jun 1966.

\bibitem[BW15]{Bridson2015b}
Martin~R. Bridson and Henry Wilton.
\newblock The triviality problem for profinite completions.
\newblock {\em Inventiones mathematicae}, 202(2):839--874, feb 2015.

\bibitem[Cei67]{Ceitin1967}
Gregory~S. Ceitin.
\newblock Algorithmic operators in constructive metric spaces.
\newblock {\em Trudy Matematicheskogo Instituta Imeni V. A. Steklova}, pages
  1--80, 1967.

\bibitem[CG05]{Champetier2005}
Christophe Champetier and Vincent Guirardel.
\newblock Limit groups as limits of free groups.
\newblock {\em Israel Journal of Mathematics}, 146(1):1--75, dec 2005.

\bibitem[Dar15]{Darbinyan2015}
Arman Darbinyan.
\newblock Group embeddings with algorithmic properties.
\newblock {\em Communications in Algebra}, 43(11):4923--4935, July 2015.

\bibitem[dCGP07]{Cornulier2007}
Yves de~Cornulier, Luc Guyot, and Wolfgang Pitsch.
\newblock On the isolated points in the space of groups.
\newblock {\em Journal of Algebra}, 307(1):254--277, jan 2007.

\bibitem[Deh11]{Dehn1911}
Max Dehn.
\newblock {\"U}ber unendliche diskontinuierliche gruppen.
\newblock {\em Mathematische Annalen}, 71(1):116--144, mar 1911.

\bibitem[Deh87]{Dehn1987}
Max Dehn.
\newblock On infinite discontinuous groups.
\newblock In {\em Papers on Group Theory and Topology}, pages 133--178.
  Springer New York, 1987.

\bibitem[DF19]{Detinko2019}
Alla~S. Detinko and Dane~L. Flannery.
\newblock Linear groups and computation.
\newblock {\em Expositiones Mathematicae}, 37(4):454--484, dec 2019.

\bibitem[DG08]{Dahmani2008}
Fran{\c{c}}ois Dahmani and Daniel Groves.
\newblock The isomorphism problem for toral relatively hyperbolic groups.
\newblock {\em Publications math{\'{e}}matiques de l'{IH}{\'{E}}S},
  107(1):211--290, June 2008.

\bibitem[DM08]{Downey2008}
Rod Downey and Antonio Montalb{\'a}n.
\newblock The isomorphism problem for torsion-free abelian groups is analytic
  complete.
\newblock {\em Journal of Algebra}, 320(6):2291--2300, September 2008.

\bibitem[EPC{\etalchar{+}}92]{Epstein1992}
David B.~A. Epstein, M.~S. Paterson, J.~W. Cannon, D.~F. Holt, S.~V. Levy, and
  W.~P. Thurston.
\newblock {\em Word Processing in Groups}.
\newblock A. K. Peters, Ltd., November 1992.

\bibitem[GMW12]{Groves2012}
Daniel Groves, Jason~Fox Manning, and Henry Wilton.
\newblock Recognizing geometric 3-manifold groups using the word problem.
\newblock October 2012.

\bibitem[Gru78]{Grunewald1978}
Fritz~J. Grunewald.
\newblock On some groups which cannot be finitely presented.
\newblock {\em Journal of the London Mathematical Society}, s2-17(3):427--436,
  June 1978.

\bibitem[GW09]{Groves2009}
Daniel Groves and Henry Wilton.
\newblock Enumerating limit groups.
\newblock {\em Groups, Geometry, and Dynamics}, pages 389--399, 2009.

\bibitem[Har11]{HARTUNG2011}
Ren{\'{e}} Hartung.
\newblock Coset enumeration for certain infinitely presented groups.
\newblock {\em International Journal of Algebra and Computation},
  21(08):1369--1380, dec 2011.

\bibitem[Her02]{Hertling2002}
Peter Hertling.
\newblock A {Banach-Mazur} computable but not {Markov} computable function on
  the computable real numbers.
\newblock In {\em Automata, Languages and Programming}, volume 132, pages
  962--972. Springer Berlin Heidelberg, mar 2002.

\bibitem[Hou07]{Houcine2007}
Abderezak~Ould Houcine.
\newblock Limit groups of equationally noetherian groups.
\newblock In {\em Geometric Group Theory}, pages 103--119. Birkhäuser Basel,
  2007.

\bibitem[HR16]{Hoyrup2016}
Mathieu Hoyrup and Crist{\'{o}}bal Rojas.
\newblock On the information carried by programs about the objects they
  compute.
\newblock {\em Theory of Computing Systems}, 61(4):1214--1236, dec 2016.

\bibitem[KMS17]{Kharlampovich2017}
Olga~G. Kharlampovich, Alexei Myasnikov, and Mark~V. Sapir.
\newblock Algorithmically complex residually finite groups.
\newblock {\em Bulletin of Mathematical Sciences}, 7(2):309--352, mar 2017.

\bibitem[Lem97]{Lempp1997}
Steffen Lempp.
\newblock The computational complexity of torsion-freeness of finitely
  presented groups.
\newblock {\em Bulletin of the Australian Mathematical Society},
  56(2):273--277, October 1997.

\bibitem[Loc81]{Lockhart1981}
Jody Lockhart.
\newblock Decision problems in classes of group presentations with uniformly
  solvable word problem.
\newblock {\em Archiv der Mathematik}, 37(1):1--6, dec 1981.

\bibitem[Mak85]{Makanin1985}
Gennadii~S. Makanin.
\newblock Decidability of the universal and positive theories of a free group.
\newblock {\em Mathematics of the {USSR}-Izvestiya}, 25(1):75--88, feb 1985.

\bibitem[Mal61]{Malcev1961}
Anatolii~I. Maltsev.
\newblock Constructive algebras {I}. ({R}ussian).
\newblock {\em Uspehi Mat. Nauk}, 1961.

\bibitem[Mal71]{Malcev1971}
Anatolii~I. Maltsev.
\newblock {\em The metamathematics of algebraic systems, collected papers:
  1936-1967}.
\newblock North-Holland Pub. Co, Amsterdam, 1971.

\bibitem[Man82]{Mann1982}
Avinoam Mann.
\newblock A note on recursively presented and co-recursively presented groups.
\newblock {\em Bulletin of the London Mathematical Society}, 14(2):112--118,
  mar 1982.

\bibitem[McK43]{McKinsey1943}
J.~C.~C. McKinsey.
\newblock The decision problem for some classes of sentences without
  quantifiers.
\newblock {\em Journal of Symbolic Logic}, 8(3):61--76, sep 1943.

\bibitem[Mel14]{MELNIKOV2014}
Alexander~G. Melnikov.
\newblock Computable abelian groups.
\newblock {\em The Bulletin of Symbolic Logic}, 20(3):315--356, September 2014.

\bibitem[Mih68]{Mihailova1968}
K.~A. Mihailova.
\newblock The occurrence problem for free products of groups.
\newblock {\em Mathematics of the USSR-Sbornik}, 4(2):181--190, February 1968.

\bibitem[{Mil}72]{MillerIII1972}
Charles~F. {Miller III}.
\newblock {\em On Group-Theoretic Decision Problems and Their Classification.
  (AM-68)}.
\newblock Princeton University Press, March 1972.

\bibitem[{Mil}92]{MillerIII1992}
Charles~F. {Miller III}.
\newblock Decision problems for groups {\textemdash} survey and reflections.
\newblock In {\em Mathematical Sciences Research Institute Publications}, pages
  1--59. Springer New York, 1992.

\bibitem[Mos64]{Moschovakis1964}
Yiannis Moschovakis.
\newblock Recursive metric spaces.
\newblock {\em Fundamenta Mathematicae}, 55(3):215--238, 1964.

\bibitem[Neu73]{Neumann1973}
Bernhard~H. Neumann.
\newblock The isomorphism problem for algebraically closed groups.
\newblock In {\em Word Problems - Decision Problems and the Burnside Problem in
  Group Theory}, pages 553--562. Elsevier, 1973.

\bibitem[OH07]{OuldHoucine2007}
Abderezak Ould~Houcine.
\newblock Embeddings in finitely presented groups which preserve the center.
\newblock {\em Journal of Algebra}, 307(1):1--23, January 2007.

\bibitem[Oza14]{Ozawa2014}
Narutaka Ozawa.
\newblock Noncommutative real algebraic geometry of {K}azhdan's {P}roperty
  {(T)}.
\newblock {\em Journal of the Institute of Mathematics of Jussieu},
  15(1):85--90, jul 2014.

\bibitem[Rau21]{Rauzy2021b}
Emmanuel Rauzy.
\newblock Computable analysis on the space of marked groups.
\newblock {\em arXiv:2111.01179}, 2021.

\bibitem[Rog87]{Rogers1987}
Hartley Rogers.
\newblock {\em Theory of Recursive Functions and Effective Computability}.
\newblock MIT Press, April 1987.

\bibitem[Sap11]{Sapir2011}
Mark~V. Sapir.
\newblock Asymptotic invariants, complexity of groups and related problems.
\newblock {\em Bulletin of Mathematical Sciences}, 1(2):277--364, mar 2011.

\bibitem[Sel01]{Sela2001}
Zlil Sela.
\newblock Diophantine geometry over groups {I}: {M}akanin-{R}azborov diagrams.
\newblock {\em Publications math{\'{e}}matiques de l'{IH}{\'{E}}S},
  93(1):31--106, sep 2001.

\bibitem[Spr98]{Spreen1998}
Dieter Spreen.
\newblock On effective topological spaces.
\newblock {\em The Journal of Symbolic Logic}, 63(1):185--221, 1998.

\bibitem[Spr16]{SPREEN2016}
Dieter Spreen.
\newblock Some results related to the continuity problem.
\newblock {\em Mathematical Structures in Computer Science}, 27(8):1601--1624,
  jun 2016.

\bibitem[Wil08]{Wilton2008}
Henry Wilton.
\newblock Hall's theorem for limit groups.
\newblock {\em Geometric and Functional Analysis}, 18(1):271--303, March 2008.

\end{thebibliography}

\end{document}